\documentclass[conf]{new-aiaa}
\usepackage[utf8]{inputenc}

\usepackage{graphicx}
\usepackage{amsmath}
\usepackage[version=4]{mhchem}
\usepackage{siunitx}
\usepackage{longtable,tabularx}
\usepackage{placeins}
\usepackage{multicol}
\usepackage{tcolorbox}
\usepackage{subcaption}
\title{Air-taxi transition trajectory optimization with\\ physics-based models}

\author{Nicholas C. Orndorff\footnote{Graduate Student, Department of Mechanical and Aerospace Engineering, AIAA Student Member.}, Luca Scotzniovsky\footnote{Ph.D Student, Department of Mechanical and Aerospace Engineering, AIAA Student Member.}, and Darshan Sarojini\footnote{Postdoctoral Scholar, Department of Mechanical and Aerospace Engineering, AIAA Member}}
\affil{University of California San Diego, La Jolla, CA, 92092}

\author{Hyunjune Gill\footnote{Ph.D Student, Department of Mechanical and Aerospace Engineering, AIAA Student Member} and Seongkyu Lee\footnote{Associate Professor, Department of Mechanical and Aerospace Engineering, AIAA Associate Fellow}}
\affil{University of California, 1 Shields Avenue, Davis, CA 95616}

\author{Zeyu Cheng\footnote{Ph.D Student, Department of Electrical and Computer Engineering}, Shuofeng Zhao\footnote{Postdoctoral Scholar, Department of Electrical and Computer Engineering}, and Chunting Mi\footnote{Professor and Chair, Department of Electrical and Computer Engineering}}
\affil{San Diego State University, San Diego, CA, 92182}

\author{John T. Hwang\footnote{Assistant Professor, Department of Mechanical and Aerospace Engineering, AIAA Member.}}
\affil{University of California San Diego, La Jolla, CA, 92092}

\begin{document}

\maketitle

\begin{abstract}
To fulfill the vision for large-scale urban air mobility, air-taxi concepts must be carefully designed and optimized for their intended mission. Proposed air-taxi missions contain dynamic segments that are dominated by nonlinear dynamics. One such segment is the transition to and from hover and cruise that occurs at the start and end of the mission. Because this transition involves low-altitude and high-power flight, analyzing transition trajectories is critical for safe and economical urban air mobility. Optimization of the transition maneuver requires an optimal control approach that characterizes the trajectories of the system states through time. In this paper we solve this optimal control problem for air-taxi transition within a large-scale design-optimization framework. This framework allows us to include five physics-based models that describe flight dynamics, rotor aerodynamics, wing aerodynamics, motor performance, and acoustics with which we create a low-fidelity model of NASA's Lift-plus-Cruise air-taxi concept. We use this optimization problem formulation to compute transition trajectories that minimize time or minimize energy. Our results show that the Lift-plus-Cruise aircraft completes a minimum-energy transition in 80s with an energy expenditure of 13.3MJ and a minimum-time transition in 28s with an energy expenditure of 16.4MJ. We find that these trajectories contain large pitch angles and high sound pressure levels which are both undesirable for practical urban air mobility. Consequently, we explore trajectories that include pitch angle and acoustic constraints, and find that minimum time trajectories are significantly more affected by these constraints than minimum energy trajectories.
\end{abstract}

\newpage
\section{Nomenclature}
\begin{multicols}{2}
\begin{center}
\begin{spacing}{0.85}

\begin{tabular}{@{}l @{\quad=\quad} l@{}}
$v$ & velocity (m/s)\\ 
$\gamma$ & flight path angle (rad)\\ 
$h$ & altitude (m)\\ 
$x$ & distance (m)\\ 
$\alpha$ & angle of attack (rad)\\
$\theta$ & pitch angle (rad)\\
$m$ & mass (kg)\\ 
$g$ & acceleration due to gravity ($\text{m/s}^{2}$)\\ 
$T_C$, $T_L$ & cruise thrust and lift thrust (N)\\ 
$L$, $D$ & lift, drag (N)\\ 
$C_L$, $C_D$ & lift/drag coefficients\\ 
$\lambda_0$, $\lambda_s$, $\lambda_c$ & rotor inflow states\\ 
\end{tabular}

\begin{tabular}{@{}l @{\quad=\quad} l@{}}
$k$ & number of lift rotors\\ 
$C_T$, $C_P$ & thrust/power coefficients \\
$P$, $Q$ & power (W) and torque (N-m)\\
$\rho$ & density ($\text{kg/m}^3$)\\
$d_C$, $d_L$ & cruise and lift rotor diameter (m)\\
$n$ & rotor speed (revolutions/s)\\
$V_{tip}$ & rotor tip speed (m/s)\\
$\sigma$ & blade solidity \\
$A_b$, $A_d$ & rotor blade and disk area ($\text{m}^2$)\\
$SPL$ & sound pressure level (db)\\
$\theta_0$ & rotor disk elevation angle (rad)\\
$s_0$ & observer distance (m)\\
\end{tabular}

\end{spacing}
\end{center}
\end{multicols}

\section{Introduction}
\lettrine{E}{lectric} air-taxi concepts are rapidly being developed to meet the vision for large-scale urban air mobility (UAM). To fulfill this vision, aircraft designers must carefully model and analyze a wide variety of systems that span a number of engineering disciplines. As such, air-taxi design is a natural candidate for large-scale multidisciplinary design optimization (MDO). Within a large-scale MDO framework, entire aircraft can be simultaneously analyzed and optimized in the presence of many (perhaps hundreds) of variables and constraints.

Large-scale MDO was applied to system-level design of electric aircraft by Hwang and Ning~\cite{hwang2018large} who optimized NASA's X-57 aircraft for multiple, steady mission segments. 
This approach is typical of classical aircraft design but lacks an analysis of the dynamic mission segments that are prevalent in UAM mission profiles such as that proposed by Silva et al.~\cite{silva2018vtol}. The primary dynamic mission segment for air taxis is the transition between hover and cruise that occurs at the start and end of most UAM missions~\cite{chauhan2020tilt,kubo2008tail,anderson2021comparison}. This transition is characterized by highly nonlinear flight dynamics and complex interactions between the aircraft and the environment. Because of this complexity, optimal transition trajectories are not necessarily intuitive. One method of finding these trajectories is to implement a trajectory optimization algorithm within an MDO framework~\cite{betts1998survey}.

Trajectory optimization problems are typically solved with one of three techniques: dynamic programming, indirect methods, and direct methods~\cite{kelly2017}. Dynamic programming relies upon a discretized state space which it explores in its entirety~\cite{sutton2018reinforcement}. Because of this, dynamic programming is best suited to small state-spaces and discrete dynamics. Complex models with continuous dynamics are often solved by indirect or direct methods. Indirect methods rely upon the analytical derivation of the necessary conditions for optimality~\cite{kelly2017}. These conditions depend upon the system model and constraints, and must be re-derived whenever the problem is modified. Direct methods start by discretizing the system and then applying optimization. Generally, the discretized system is transcribed as a nonlinear programming problem (NLP) that can be solved using gradient-based optimization~\cite{hargraves1987direct,bryson1962steepest}. With direct transcription there is no need to derive the necessary conditions for optimality~\cite{betts1998survey}. Because of this, direct transcription is well suited to complex systems where constraints and models rapidly evolve (e.g., the aircraft design process).

The conversion from an optimal control problem to a nonlinear program requires both discretization and simulation. Generally, continuous systems are discretized by introducing collocation points that represent a finite position in the state space where controls and constraints can be defined. The state values at the collocation points are defined through simulation. Simulating trajectories by explicitly propagating states across the collocation points is known as "shooting" and typically makes use of standard ODE solution methods~\cite{betts1998survey}.
Hwang and Munster~\cite{hwang2018solution} applied the general linear methods (GLM) approach in the implementation of an ODE solver within an  MDO framework. The GLM method simplifies the implementation of collocation within the NLP framework by allowing rapid implementation of multiple integration methods. Derivatives of the resulting ODE representation must be provided in order to apply gradient-based optimization. These can be computed using the modular analysis and unified derivatives (MAUD) approach~\cite{martins2013review,hwang2018computational}.

Direct transcription is commonly used for aircraft trajectory optimization in the presence of multiple physics-based models and constraints. Betts and Cramer~\cite{betts1995application} demonstrated direct transcription for a commercial aircraft trajectory with a large number of realistic constraints derived from federal airspace regulations. Chauhan and Martins~\cite{chauhan2020tilt} used direct transcription to find optimal trajectories for the Airbus Vahana tilt-wing aircraft. They evaluated the effects of varying acceleration and stall constraints in order to explore optimal transition trajectories. Hargraves et al.~\cite{hargraves1981numerical} demonstrated the versatility and robustness of direct methods by finding optimal trajectories for a variety of aerospace applications. These included minimum-time climb and minimum-fuel missions, as well as rocket ascent and evasive maneuvers. Recent work by Hendricks et al.~\cite{hendricks2019multidisciplinary} and Jasa et al.~\cite{jasa2020large} exemplified the use of the MAUD and GLM methods with applications for air taxis and supersonic aircraft respectively.


One advantage of direct transcription is the ability to rapidly evaluate sub-models. Anderson et al.~\cite{anderson2021comparison} evaluated the effect of several aerodynamic sub-models on transition trajectories for a small bi-wing tailsitter aircraft. They specifically evaluated rotor-wing interactions, and found that higher fidelity models predicted more conservative trajectories. Falck et al.~\cite{falck2018multidisciplinary} implemented an acoustic sub-model with which they found noise-limited trajectories for a quadrotor air taxi. From these papers we deduce that model selection has a significant effect on the validity of the system trajectories.

In this paper we use direct transcription to find optimal transition trajectories for NASA's Lift-plus-Cruise air-taxi concept~\cite{silva2018vtol}. The Lift-plus-Cruise aircraft exemplifies the multidisciplinary nature of UAM aircraft design by combining a fixed wing with lift and cruise rotors (Fig.~1). We create five sub-discipline models to accurately describe the aircraft: flight dynamics, rotors aerodynamics, wing aerodynamics, motor performance, and aeroacoustics. We combine these models within a nonlinear programming problem in which we explicitly propagate the states through time to obtain their trajectories. This formulation allows us to rapidly evaluate different constraints and objectives. Specifically, we explore constraints on pitch angle and sound pressure level with respect to two objectives: minimum-energy and minimum-time transition trajectories.

\FloatBarrier
\begin{figure}[ht]
\centering
\includegraphics[width=0.8\textwidth]{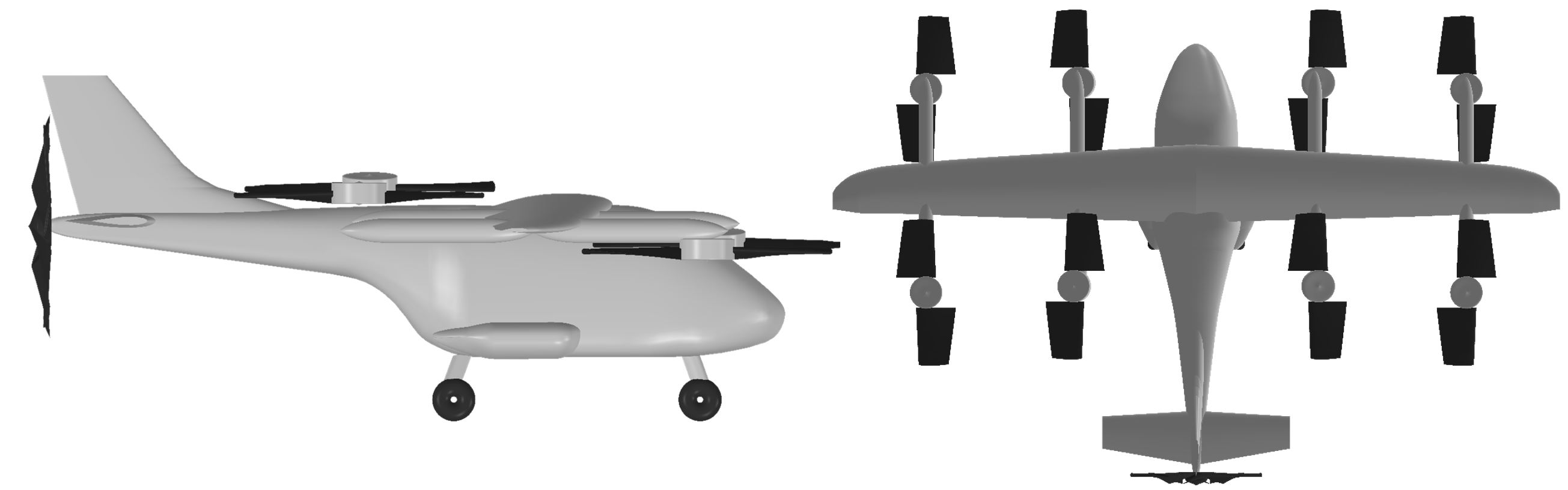}
\caption{NASA's Lift-plus-Cruise air-taxi concept~\cite{silva2018vtol}.}
\end{figure}
\FloatBarrier

\section{Methodology}
\subsection{Flight Dynamics}
Aircraft trajectories are simulated with a three-degree-of-freedom system of ODEs. The following four equations allow for the explicit propagation of velocity $(v)$, flight path angle $(\gamma)$, altitude $(h)$, and horizontal position $(x)$ with respect to time:

\begin{equation}\label{dv}
\dot{v} = \frac{T_C}{m}\cos{\alpha} + \frac{T_L}{m}\sin{\alpha} - \frac{D}{m} - g\sin{\gamma}
\end{equation}

\begin{equation}\label{dgamma}
\dot{\gamma} = \frac{T_C}{mv}\sin{\alpha} + \frac{T_L}{mv}\cos{\alpha} + \frac{L}{mv} - \frac{g\cos{\gamma}}{v}
\end{equation}

\begin{equation}\label{dh}
\dot{h} = v\sin{\gamma}
\end{equation}

\begin{equation}\label{dx}
\dot{x} = v\cos{\gamma}.
\end{equation}

Figure 2 depicts these four aircraft states and their corresponding coordinate systems. The inputs to the system of ODEs are the lift and cruise rotor thrust ($T_C$, $T_L$), the angle of attack ($\alpha$), lift ($L$), and drag ($D$). In previous work~\cite{orndorff2022investigation} we found that (a) the difference in thrust and power across the eight lifting rotors is negligible, and (b) the difference in inflow velocities due to aircraft rotation is also negligible. Consequently, we reduce the computational complexity of the model by analyzing a single lift rotor and multiplying the thrust and power by $k$ lift rotors to get the total lift rotor thrust and the total lift rotor power (e.g., $T_L=kT$).

\begin{figure}[ht]
\centering
\includegraphics[width=.55\textwidth]{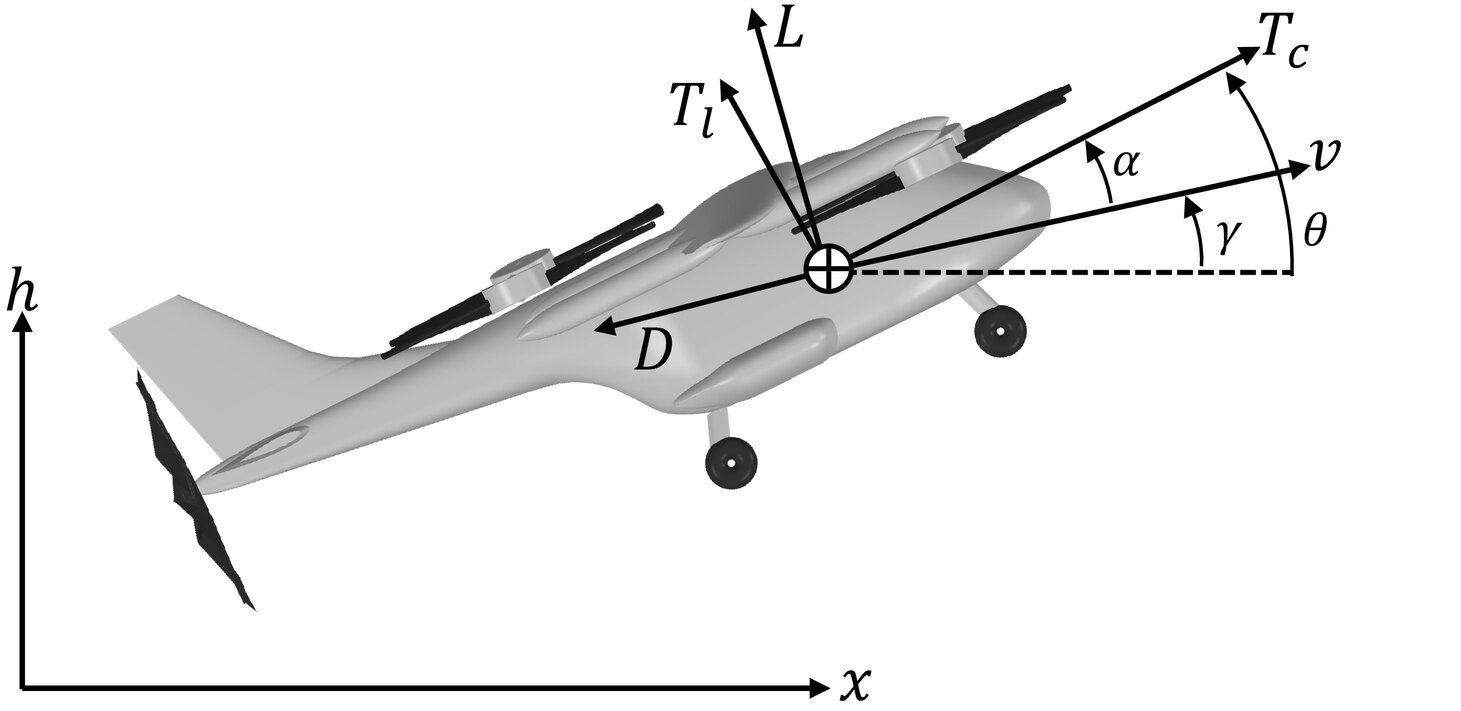}
\caption{Coordinate system for the 3-DOF flight dynamics model.}
\end{figure}

\subsection{Rotor Aerodynamics}
Air-taxi transition trajectories inherently contain large segments of edgewise flight. We account for this by implementing the dynamic inflow rotor model developed by Pitt and Peters~\cite{pitt1980rotor,pitt1981theoretical}. The Pitt-Peters model accounts for edgewise flight with three inflow states ($\lambda_0$, $\lambda_s$ and $\lambda_c$) which are the uniform, side-to-side, and fore-to-aft flow variations respectively. The rotor induced downwash $(\lambda_i)$ is calculated as a function of the dimensionless radial distance from the rotor hub $r^{*}$ and the azimuth angle of the rotor disk $\psi$:

\begin{equation}
    \lambda_i = \lambda_0 + \lambda_c r^{*}\cos\psi +\lambda_s r^{*}\sin\psi.
\end{equation}

In the following system of equations, these dynamic inflow parameters are combined with coefficients calculated with blade element theory ($C_T$, $C_{M_{y}}$, $C_{M_{x}}$), the derivative matrix $L$, and the mass matrix $M$:

\begin{equation}
    \label{Pitt-Peters differential equation}
    M \begin{bmatrix}
       \dot{\lambda_0} \\
       \dot{\lambda_c} \\
       \dot{\lambda_s} \\
    \end{bmatrix}
    + L^{-1} \begin{bmatrix}
        \lambda_0 \\
        \lambda_c \\
        \lambda_s \\
    \end{bmatrix}
    = \begin{bmatrix}
        C_T \\
        -C_{M_{y}}\\
        C_{M_{x}} \\       
    \end{bmatrix}.
\end{equation}

To simplify the implementation of this rotor model, we use surrogate models for the thrust and power coefficients $(C_T \text{ and } C_P)$. The surrogate models are created with regularized minimal-energy tensor-product splines (RMTS)~\cite{hwang2018fast} using the surrogate modeling toolbox (SMT) developed by Bouhlel et al.~\cite{SMT2019}. Figure 3 shows the surrogate models created for the cruise rotor.

\begin{figure}[ht]
\centering
\includegraphics[width=1\textwidth]{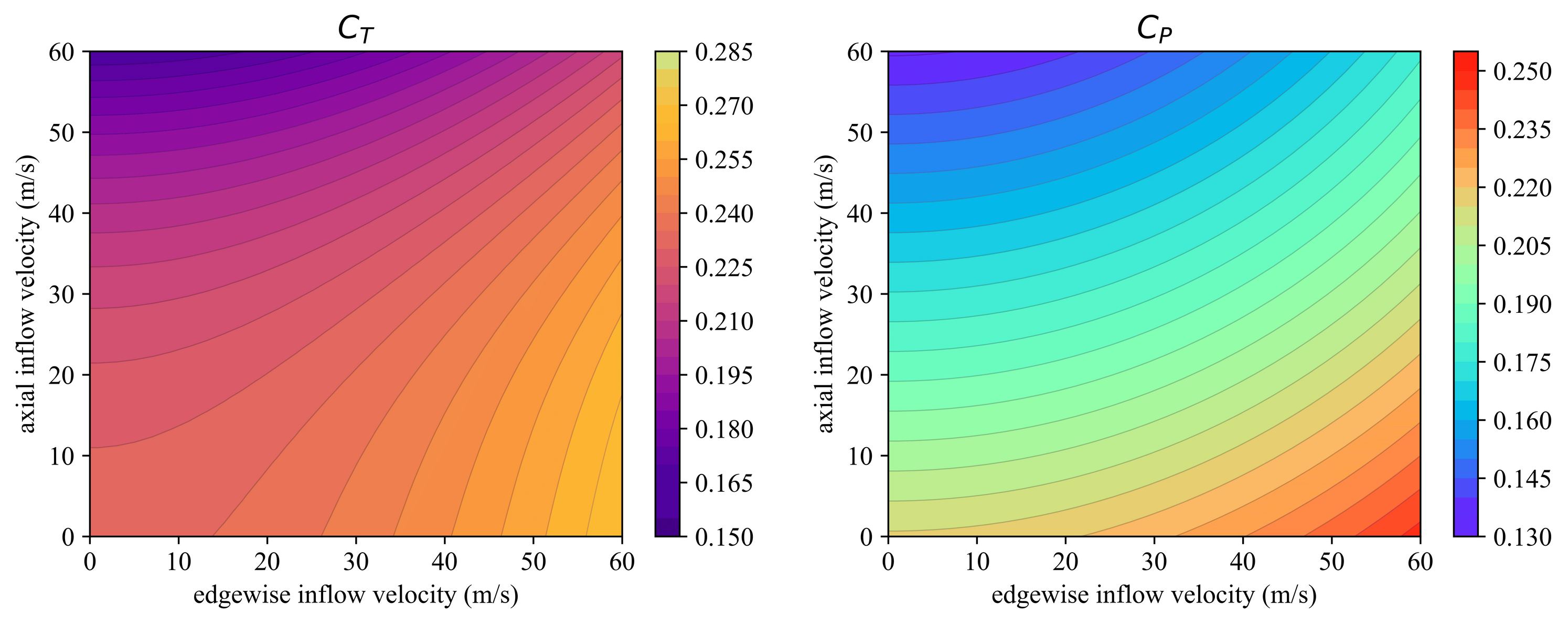}
\caption{Cruise rotor surrogate models for thrust and power coefficients as a function of axial and edgewise free-stream velocity.}
\end{figure}

Rotor thrust $(T)$, torque $(Q)$, and power $(P)$ are calculated with these surrogate models as functions of the rotor diameter $(d)$ (see $d_C$ and $d_L$ in Table 1), rotations per second $(n)$, and atmospheric density $(\rho)$:

\begin{equation}
T = \rho n^2 d^4 C_T
\end{equation}
\begin{equation}
P = \rho n^3 d^5 C_P
\end{equation}
\begin{equation}
Q= \rho n^2 d^5 C_Q.
\end{equation}

The torque coefficient $C_Q$ is found as a function of the power coefficient: $C_Q=C_P/2\pi$. To compute the total aircraft power we sum the individual rotor powers (with $k$ lifting rotors): $P = P_{\text{Cruise}} + kP_{\text{Lift}}$. The total energy $E$ consumed by the aircraft is expressed as an additional ordinary differential equation $\dot{E}=P$. The total energy used by the aircraft for a given trajectory is the final entry in the resulting discrete-time energy vector.

\subsection{Wing Aerodynamics}
Wing aerodynamics are analyzed using standard lift and drag approximations. Lift and drag coefficients are found in a surrogate model that spans a wide range of angle of attack values ($-90^{\circ}$ to $90^{\circ}$) (Fig. 4). These models are constructed with radial basis functions (RBF) using the SMT toolbox~\cite{SMT2019} based in part on data from XFOIL~\cite{drela1989xfoil} and extrapolated to high $\alpha$ using data from Selig~\cite{selig2014real}. In the following equations, $s$ is the wing area, and $\rho$ is the atmospheric density computed, as above, from a standard atmospheric model:

\begin{equation}
L = \frac{1}{2}\rho v^2 s C_L
\end{equation}
\begin{equation}
D = \frac{1}{2}\rho v^2 s C_D.
\end{equation}

Specific aircraft parameters required for the aerodynamics model can be found in Table 1. Using the aspect ratio and wing area we apply a finite wing correction by including the span efficiency ($e=0.9$) when computing the lift coefficient $C_L$ and the induced drag $C_{Di}$. The angle of attack at the wing is the sum of the wing incidence and the aircraft angle of attack ($\alpha_w=i_w + \alpha$) and the drag coefficient is found with a typical drag buildup ($C_D=C_{D0} + C_{Di}$).

\begin{figure}[ht]
\centering
\includegraphics[width=1\textwidth]{wing_model}
\caption{Aerodynamic surrogate model for lift and drag coefficients (NACA-2412).}
\end{figure}

\begin{table}[ht]
\begin{spacing}{0.85}
\caption{\label{aircraft} Lift-plus-Cruise air-taxi parameters}
\centering
\begin{tabular}{cc}
\hline \hline
Mass ($m$) & 2000kg \\
Wing area ($s$) & $19.6\text{m}^2$ \\
Aspect ratio ($AR$) & 12.13 \\
Wing incidence ($i_w$) & $2^{\circ}$ \\
Zero-lift drag coefficient ($C_{D0}$) & 0.02 \\
Cruise / lift rotor diameters ($d_C$, $d_L$) & 2.0m / 2.2m \\
Mean rotor chord ($\text{MAC}_r$) & 0.15m \\
\hline \hline
\end{tabular}
\end{spacing}
\end{table}

\subsection{Aeroacoustics}
Noise constraints on air-taxi operation present a major barrier to large-scale UAM adoption~\cite{bryson2016multidisciplinary}. Rotor noise is broadly categorized as tonal noise and broadband noise. In particular, it was found that rotor broadband noise significantly affects human hearing perception and plays a significant role in the public acceptance of UAM aircraft~\cite{Lee:2021:PAS,Li:2020:AHS,Li:2021:JAHS,Li:2022:JAHS,Li:2022:AHS,Li:2022:AIAA,Greenwood:2022:IJA}. For this paper, we only consider broadband noise due to its importance and relative simplicity. We compute broadband noise using the model developed by Schlegel, King, and Mull~\cite{schlegel1966helicopter}:


\begin{equation}\label{spl150}
SPL_{150} = 10\log_{10}{[V_{tip}^6 A_b (C_t/\sigma)^2]} - 42.9.
\end{equation}

Equation \eqref{spl150} calculates the sound pressure level (SPL) 150m below the rotor as a function of the blade tip speed ($V_{tip}$) and the blade loading ($C_t/\sigma$). The blade solidity ($\sigma$) is the ratio of blade area to disk area ($\sigma\equiv A_b/A_d$). The disk area is calculated as $A_d=\pi(d/2)^2$ and the blade area is approximated by $A_b=n_b \text{MAC}_r R$ where $n_b$ is the number of blades and $\text{MAC}_r$ is the mean rotor chord from Table 1. The lift rotors operate in edgewise flight; therefore, we include an airspeed correction term $(f(\lambda)\approx 8.3\lambda)$~\cite{Davidson:1965:RSS} as a function of the advance ratio $(\lambda=v/(nd))$ when computing the lift rotor noise:

\begin{equation}
SPL_{150} = 10\log_{10}{[V_{tip}^6 A_b (C_t/\sigma)^2]} + f(\lambda) - 42.9.
\end{equation}

We use the equation from Johnson~\cite{johnson2013rotorcraft} to compute the sound pressure level at a distance of one rotor diameter $(d)$ directly below the lift and cruise rotors:

\begin{equation}
SPL_{d} = SPL_{150} + 20\log_{10}{\left(\frac{150}{d}\right)}.
\end{equation}

The sound pressure level at an arbitrary observer location can then be predicted using the directivity relationship presented by Li and Lee~\cite{li2020prediction}. In the following equation the elevation angle from the rotor plane to the observer location is denoted by $\theta_0$ and the observer distance by $s_0$ (see Fig. 5):

\begin{equation}
SPL = |\sin{\theta_0}|^{\beta_1}SPL_d - (\beta_2 + \beta_3(1 - |\sin{\theta_0}|))\log_{10}{\bigg(\frac{|s_0|}{d}\bigg)}.
\end{equation}

The three values $\beta_1$, $\beta_2$, and $\beta_3$ are 0.0209, 18.2429, and 6.7267, respectively~\cite{li2020prediction}. Note that the maximum sound pressure level at ground level is not necessarily directly beneath the aircraft. We address this by calculating SPL within $\pm500$m of the current horizontal position $x_i$, $\forall t_i \in \vec{t}$ (Fig. 5). The ground-level sound pressure for the single cruise rotor and $k$ lift rotors are summed in Eq.~\eqref{splsum}. We then calculate the maximum SPL over the $x_i\pm500$m range ($\text{max}[SPL_{sum}]$), which is taken as the ground-level SPL for the current time step within the aircraft's trajectory.

\begin{equation}\label{splsum}
SPL_{sum} = 10\log_{10}{\left(10^{\frac{SPL_{d,cruise}}{10}} + k10^{\frac{SPL_{d,lift}}{10}} \right)}.
\end{equation}

\FloatBarrier
\begin{figure}[ht]
\centering
\begin{subfigure}[h]{0.4\linewidth}
\includegraphics[width=1\linewidth]{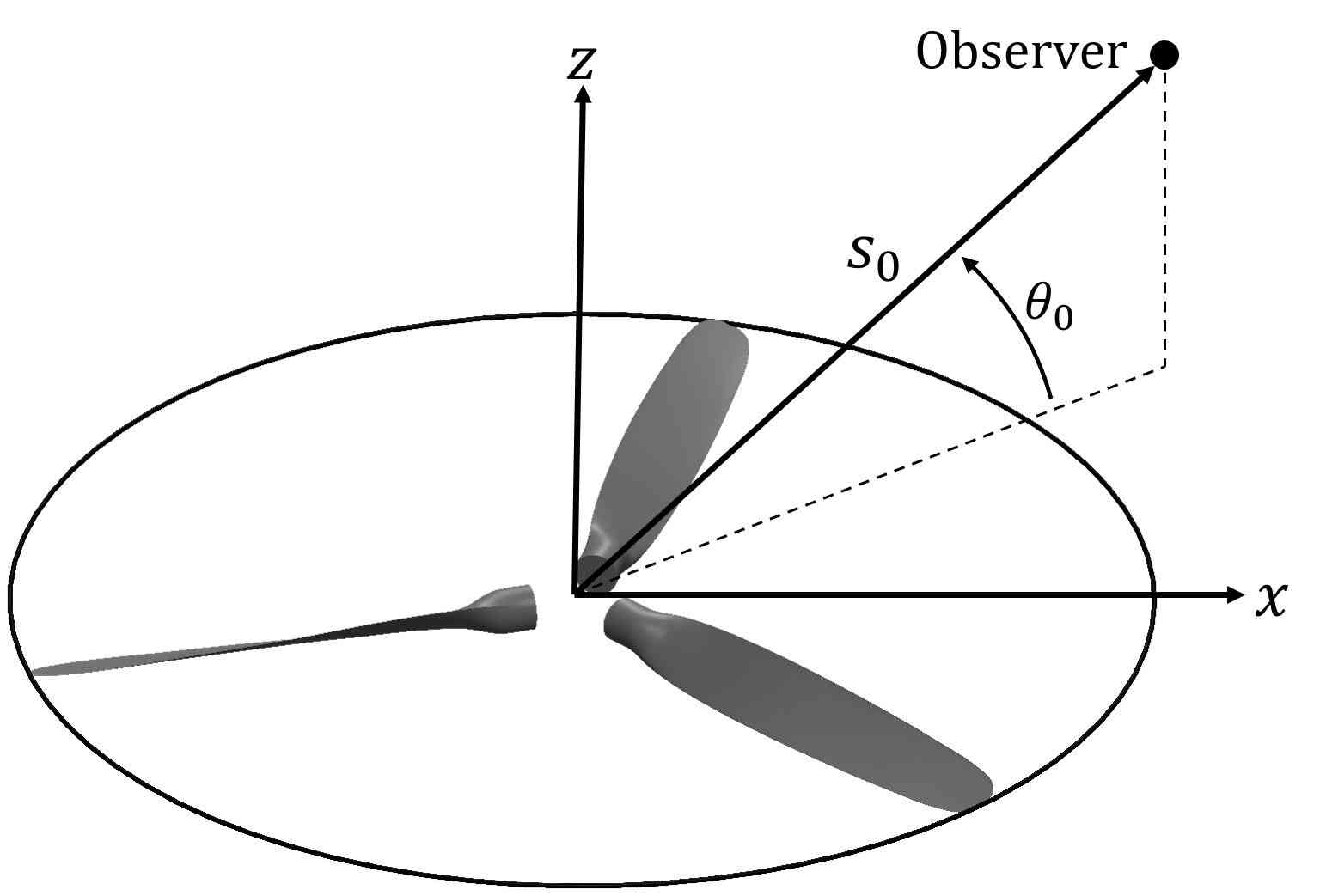}
\end{subfigure}
\begin{subfigure}[h]{0.4\linewidth}
\includegraphics[width=1\linewidth]{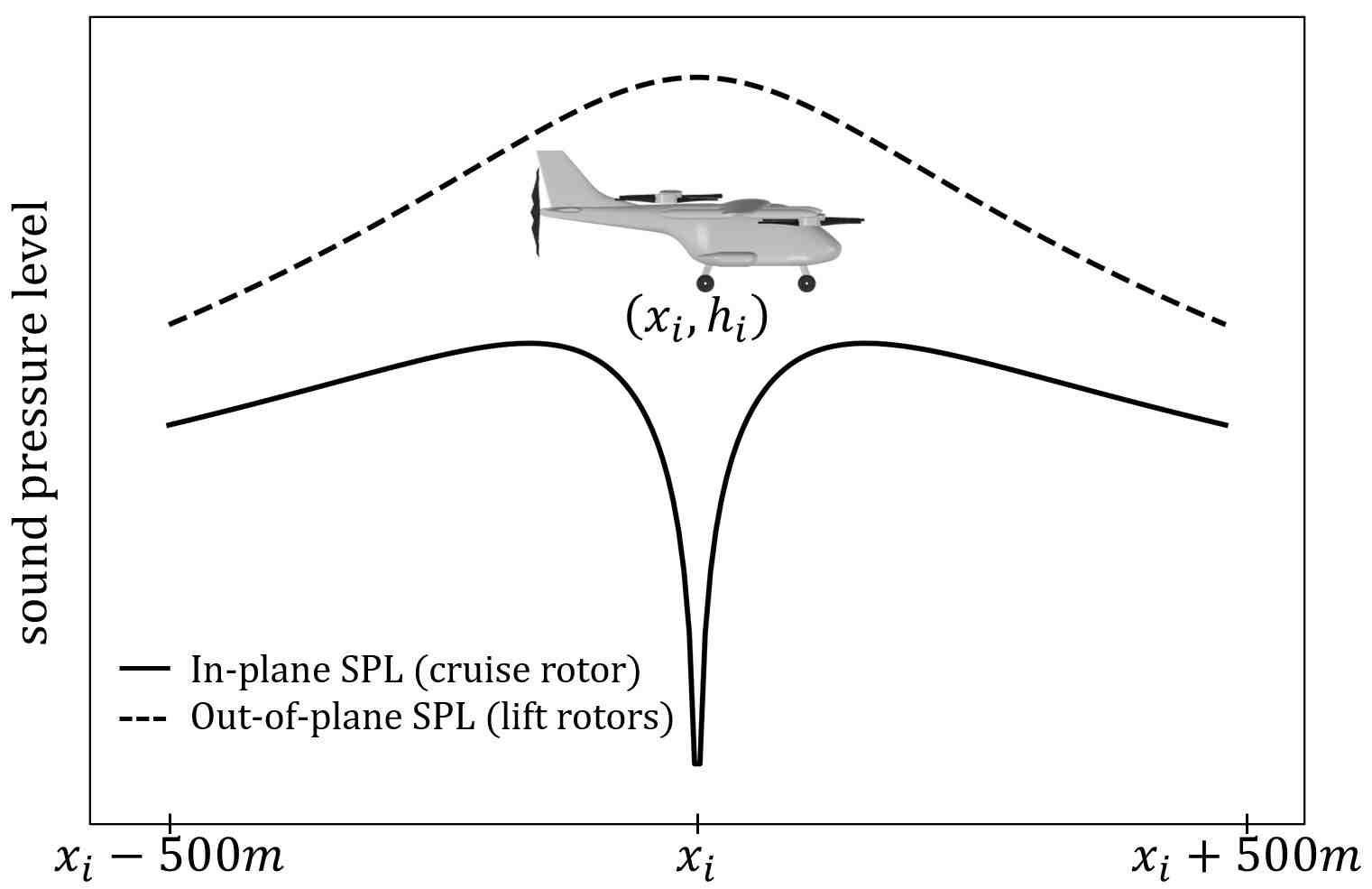}
\end{subfigure}
\caption{Rotor broadband noise diagram for an arbitrary observer location (left). At each time step the sound pressure level from each rotor is evaluated over $\pm \text{500m}$ (right).}
\end{figure}

We note that, in the current paper, tonal noise is not included in the trajectory optimization. However, tonal noise can significantly influence the overall SPL of multi-rotor UAM aircraft~\cite{jia2022computational,Jia:2022:JAHS,Sagaga:2021:VFS,Gill:2023:AIAA}. The acoustically constrained results presented here would certainly be different with the inclusion of tonal noise, which will be investigated in future work.

\subsection{Motors}

We account for the effects of motor efficiency by including a two-part model. First, we include a sizing model based on empirical relationships~\cite{htwe2019design}. This model calculates internal motor geometry as a function of two high-level parameters: stator diameter and motor length. The outputs of the sizing model are inputs to the second model, which calculates the motor efficiency with respect to load torque and motor speed by solving the following residual:

\begin{equation}
r =  \eta \tau_{EM} - \tau_{load} = 0.
\end{equation}

The motor analysis model is represented in an implicit manner in order to relate the unknown electromagnetic torque $\tau_{EM}$ and efficiency $\eta$ to the load torque $\tau_{load}$. The motor power and current are calculated with two discrete control modes: flux weakening~\cite{lu2010review} and maximum torque per ampere (MTPA)~\cite{zhao2015control}. These modes govern motor operation at the voltage limit $(V_{lim}=400)$ and below the voltage limit respectively. We can then use these models to generate motor efficiency maps (Figure 6). Using the motor efficiency map we calculate the total power as a function of the ideal rotor power $P=P_r/\eta$.

\begin{figure}[ht]
\centering
\includegraphics[width=0.55\textwidth]{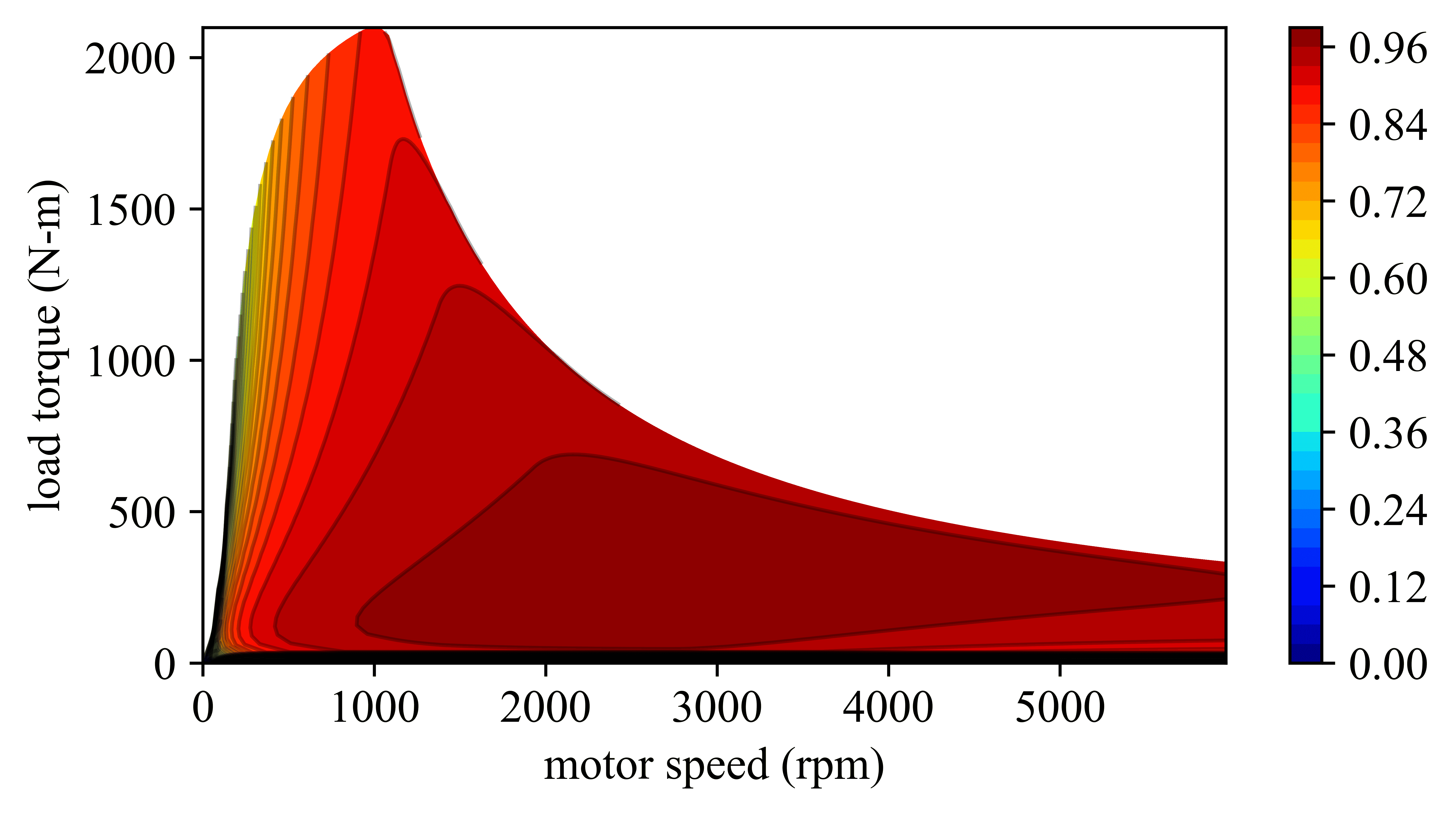}
\caption{Motor efficiency map.}
\end{figure}

\FloatBarrier
\subsection{Trajectory Optimization}
We begin by describing the trajectory optimization problem with an optimal control representation. The system dynamics are described in a general form: $\dot{\mathbf{x}}(t)=f(t,\mathbf{x}(t),\mathbf{u}(t))$ and we aim to find control profiles $\mathbf{u}(t)$ that minimize a certain objective function $J(t,\mathbf{x}(t),\mathbf{u}(t))$. The resulting state trajectories $\mathbf{x}(t)$ must satisfy certain path and terminal constraints ($c_p$ and $c_t$). The complete optimal control representation is summarized as:

\begin{equation}
\begin{array}{r l}
\underset{\mathbf{u}(t)}{\text{minimize}} & J(t,\mathbf{x}(t),\mathbf{u}(t))\\
\text{subject to} & c_p(t,\mathbf{x}(t),\mathbf{u}(t))\leq0 \\
& c_t(t_f,\mathbf{x}(t),\mathbf{u}(t))=0. \\
\end{array}
\end{equation}

We evaluate two different objectives: minimum energy and minimum time transition trajectories. Energy is calculated as the integral of total motor power $E=\int P(t,x(t),u(t))dt$ and time is represented by the integrator time step $\Delta t$ where $t_f=N\Delta t$ for $N$ time steps. Using direct transcription we reformulate the optimal control problem in a discretized form for the minimum energy objective:

\begin{equation}
\begin{array}{r l l}
\underset{\vec{u}, \Delta t}{\text{minimize}}&E(\vec{x},\vec{u})=\sum_{t=0}^{t_f} P(x_t, u_t), \forall \vec{t},\vec{x},\vec{u} \in \mathcal{R}^{N+1}& \text{(objective function)}\\
\text{subject to}&c_p(\vec{x},\vec{u})\leq0&\text{(path constraints)}\\
&c_t(\vec{x},\vec{u}) = 0&\text{(terminal constraints)}\\
\text{with}&\vec{t} = [t_i, t_i + \Delta t, t_i + 2\Delta t, \ldots, t_i + N \Delta t]&\text{(discrete time vector).} \\
\end{array}
\end{equation}

We consider three dynamic variables as the inputs $\vec{u}$: angle of attack, cruise rotor speed, and lift rotor speed. We include a single static input: the integrator time step $\Delta t$. Since the number of time steps is a fixed quantity $(N)$, $\Delta t$ is effectively a method of adjusting the total transition time $t_f$.

The aircraft trajectories are simulated with a fourth-order Runge-Kutta integration scheme using initial conditions that are derived from a steady hovering condition ($v_0=0$m/s, $\gamma_0=0^{\circ}$, $x_0=0$m, and $h_0=0$m). The transition maneuver is considered complete when several terminal constraints are satisfied ($h_f=300$m and $v_f=43$m/s). Path constraints are comprised of a variety of practical aircraft and trajectory limitations. For example, rotor power is constrained by the maximum available power, and the minimum altitude is constrained by ground collisions. All path constraints are implemented using maximum and minimum operators across state vectors in order to reduce the total number of constraints and increase the optimization speed. The complete optimization problem is summarized in Table 2.

\begin{table}[ht]
\begin{spacing}{0.85}
\caption{\label{optsum} Optimization problem formulation}
\centering
\begin{tabular}{cccc}
\hline \hline
\textsc{Objective} & Energy or time & $E$, $\Delta t$ & \textsc{Quantity} ($N$ timesteps)\\ 
\hline
\textsc{Design variables} & Lift rotor speed & $n_L$ & $N$\\ 
 & Cruise rotor speed & $n_C$ & $N$\\
 & Angle of attack & $\alpha$ & $N$\\
 & Timestep & $\Delta t$ & 1\\
 & & & Total design variables: $(3N + 1)$\\
\hline
\textsc{Constraints} & Lift power & $\max{\vec{P_L}}\leq103,652W$ & 1\\
 & Cruise power & $\max{\vec{P_C}}\leq468,300W$ & 1\\
 & Final altitude & $h_f=300\text{m}$ & 1\\
 & Minimum altitude & $\min{\vec{h}}\geq h_0$ & 1\\
 & Final velocity & $v_f=43\text{m/s}$ & 1\\
 & & & Total constraints: $5$\\
\hline
\textsc{Other constraints} & Sound pressure level & $\max{\vec{SPL}}\leq SPL_{max}$ & $N$\\
 & Pitch angle & $|\theta|\leq\theta_{max}$ & 1\\
\hline \hline
\end{tabular}
\end{spacing}
\end{table}

We solve the discretized optimization problem with gradient-based optimization and adjoint-based derivative computation. The computational expense of adjoint derivatives scales well with large numbers of design variables, and is well suited to large-scale MDO. The optimizer is SNOPT-7.7~\cite{gill2005snopt} (a sparse nonlinear optimizer), which in most cases converges to an optimal solution within one hour using a modern laptop with an Intel i7 processor and 16GB of memory. The system is discretized with $N=30$ time steps, a number that is empirically determined to provide sufficient resolution without unnecessary computation time. This results in an optimization problem with 91 design variables and between 5 and 35 constraints.

The importance of variable scaling for NLP problems cannot be understated~\cite{gills1981practical}. Certain variables such as power are large $(\mathcal{O}(10^6))$ while others can be quite small such as angle of attack $(\mathcal{O}(1^{-1}))$. This leads to a poorly conditioned problem, that struggles to converge. We address this by scaling variables before optimization. For variables that are roughly constant, we accomplish this by multiplying the variable by the inverse of the expected value. In other cases we use an affine transformation~\cite{patterson2014gpops} $\tilde{x}=vx+r$ to scale an arbitrary variable $x\in[a,b]$ to $\tilde{x}\in[-1,1]$. This uses a scale value $v = \frac{1}{b-a}$ and a shift value $r = 1-\frac{b}{b-a}$.

\section{Results}
Using the problem formulation described above, we now present optimal transition trajectories for the minimum energy and minimum time objectives. We characterize these trajectories by plotting the time histories of the system states $(v,\gamma,h,x)$. To gain further insight, we also plot the control inputs $(\alpha,n_L,n_C)$ and functions of the states (e.g., power and SPL).

To facilitate comparisons between trajectories, we use a transition efficiency parameter $\eta$ as defined by Chauhan and Martins~\cite{chauhan2020tilt}. This parameter computes the ratio of the increase in mechanical energy (kinetic plus potential energy) to the total energy consumed by the aircraft: $\eta=(0.5m\Delta v^2 + mg\Delta h)/E$. This equation is equivalent to the ratio of the specific energies (energy per unit mass), and enables future comparisons between different aircraft configurations.

We begin this section by presenting the results of the nominal case. That is, we find optimal transition trajectories resulting from the problem description in Table 2 without the additional constraints. We then apply the constraints on pitch angle and SPL. These constrained results are compared with the nominal case for both objectives, and the implications are discussed.

\subsection{Optimal Transition Trajectories for Minimum Time and Minimum Energy Objectives}

We start by plotting altitude versus horizontal displacement $(h \text{ vs. } x)$ from the point of view of an observer in an Earth-fixed frame (Fig. 7). This is instructive, as it describes the trajectory that an actual aircraft would follow if minimum energy/time trajectories were required. Adjacent data points in these figures are separated by the time step $\Delta t$.

\begin{figure}[ht]
\centering
\includegraphics[width=0.9\textwidth]{new_h_vs_x.png}
\caption{Altitude versus horizontal displacement for minimum energy and minimum time trajectories.}
\end{figure}

The minimum energy trajectory in Fig. 7 exhibits two distinct segments. First, the aircraft accelerates with little altitude variation through a horizontal distance of
 800m. Second, the aircraft climbs at a roughly steady rate ($5.5$m/s) before abruptly levelling off at the target altitude of 300m. The total horizontal distance traveled is 3,476m, which significantly eclipses the horizontal distance covered by the minimum time trajectory (1,122m). The minimum time trajectory is characterized by a bump in the trajectory around 200m of horizontal displacement. This behavior coincides with a rapid decrease in flight path angle (Fig. 8) and a slight plateau in the lift rotor speed (Fig. 9). Based on this evidence, we conclude that the bump signifies a crossing over point between lift produced primarily by the lifting rotors and lift produced by the fixed wing. Future work will include a constraint on the final flight path angle.

\FloatBarrier
\begin{figure}[ht]
\centering
\begin{subfigure}[h]{0.45\linewidth}
\includegraphics[width=1\linewidth]{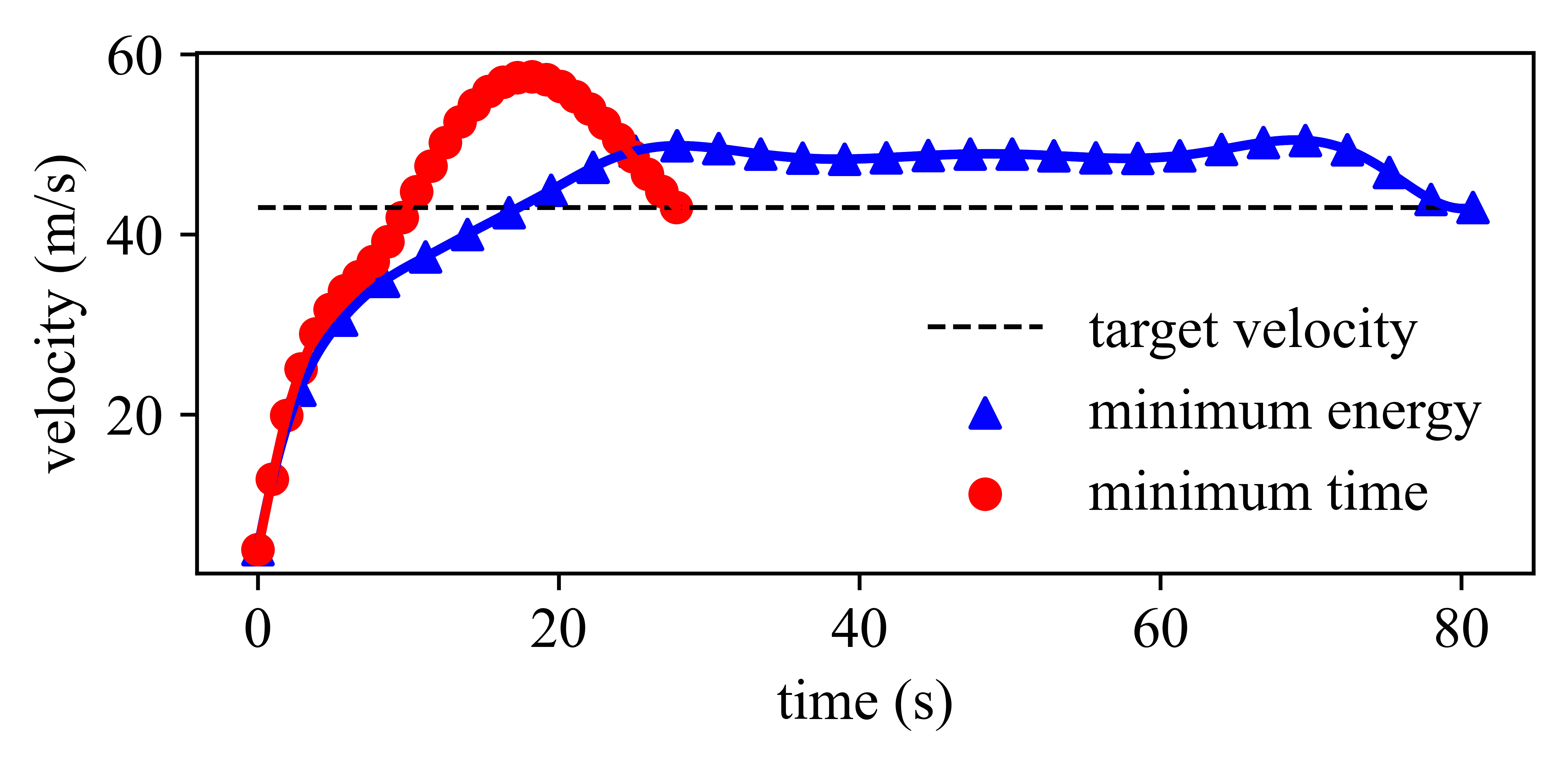}
\end{subfigure}
\begin{subfigure}[h]{0.45\linewidth}
\includegraphics[width=1\linewidth]{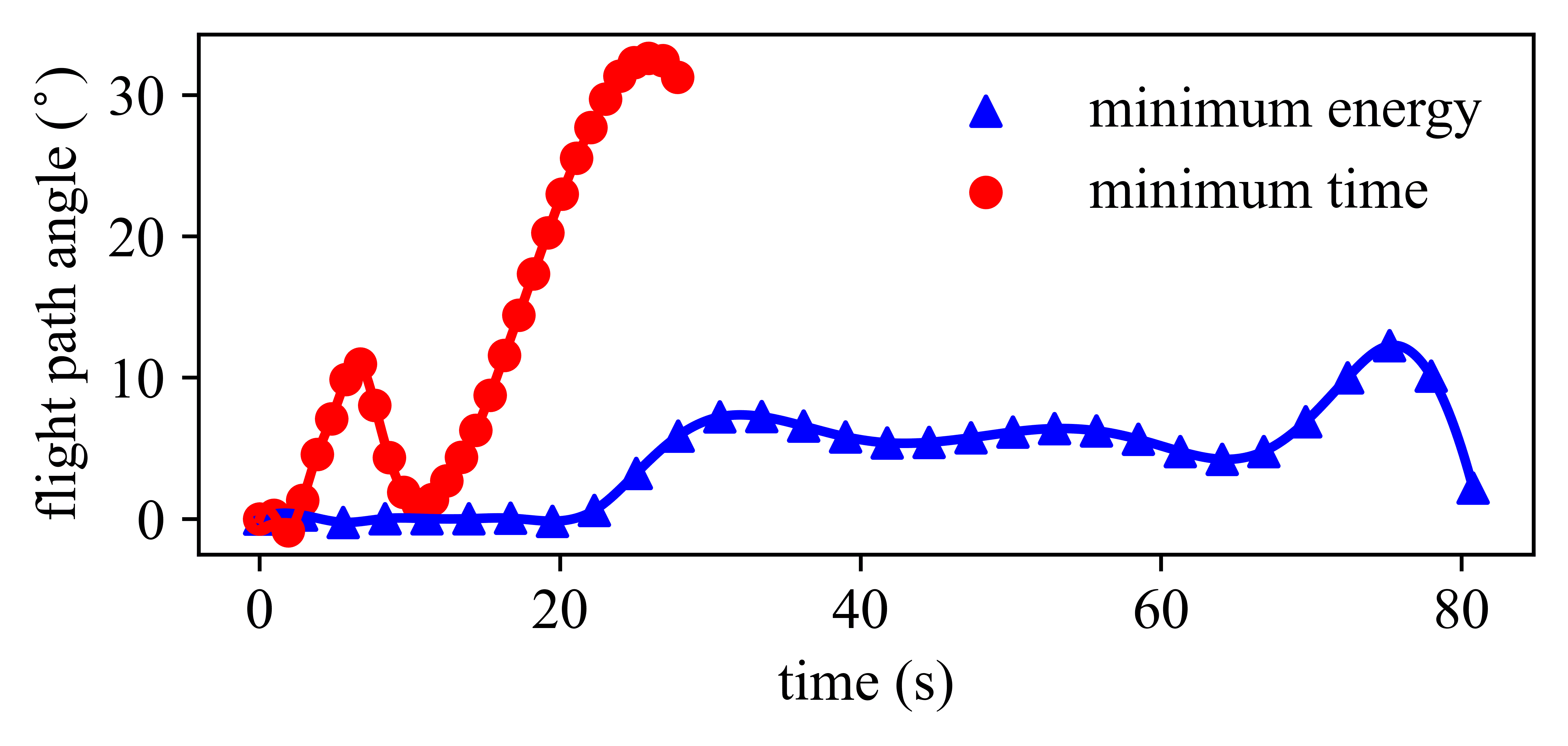}
\end{subfigure}
\begin{subfigure}[h]{0.45\linewidth}
\includegraphics[width=1\linewidth]{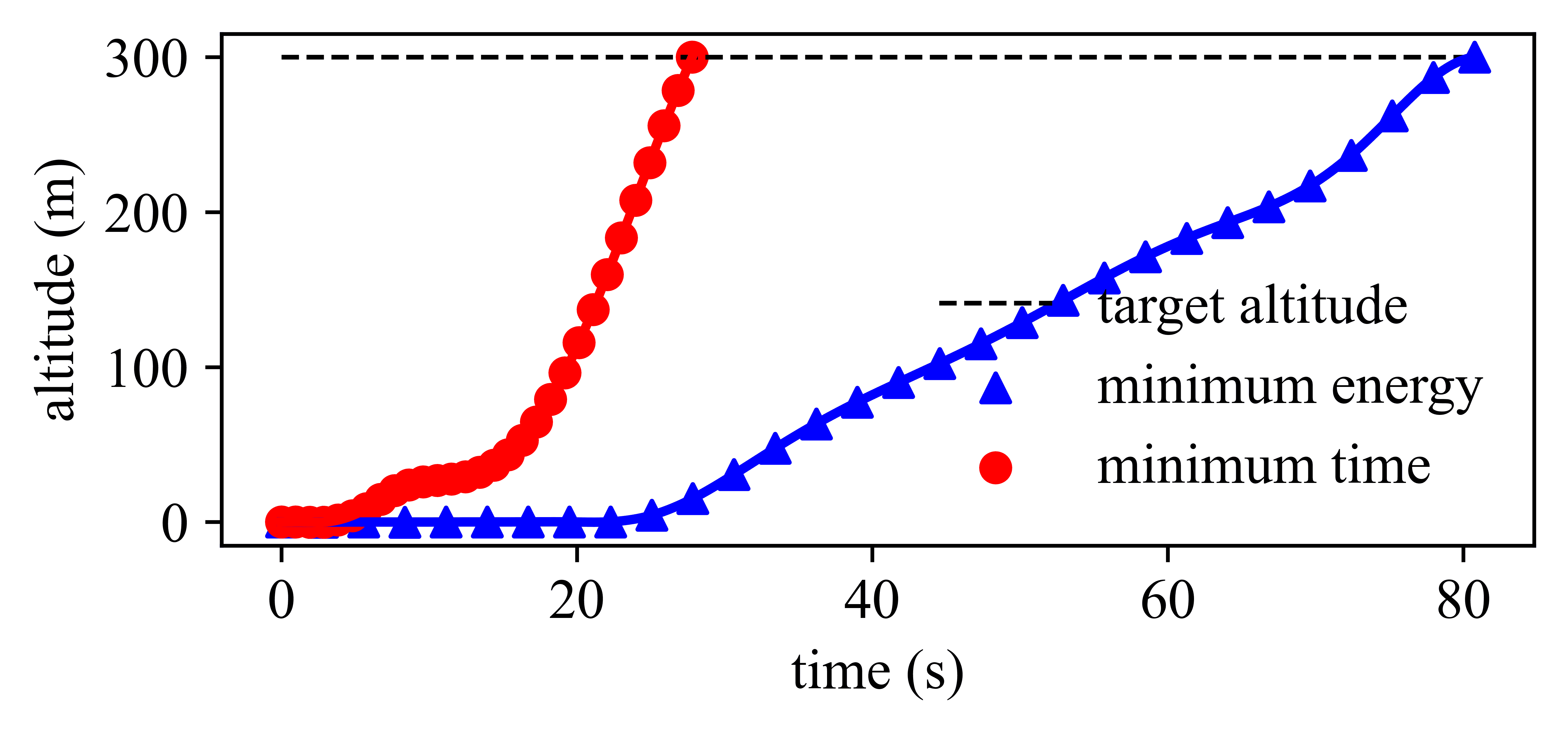}
\end{subfigure}
\begin{subfigure}[h]{0.45\linewidth}
\includegraphics[width=1\linewidth]{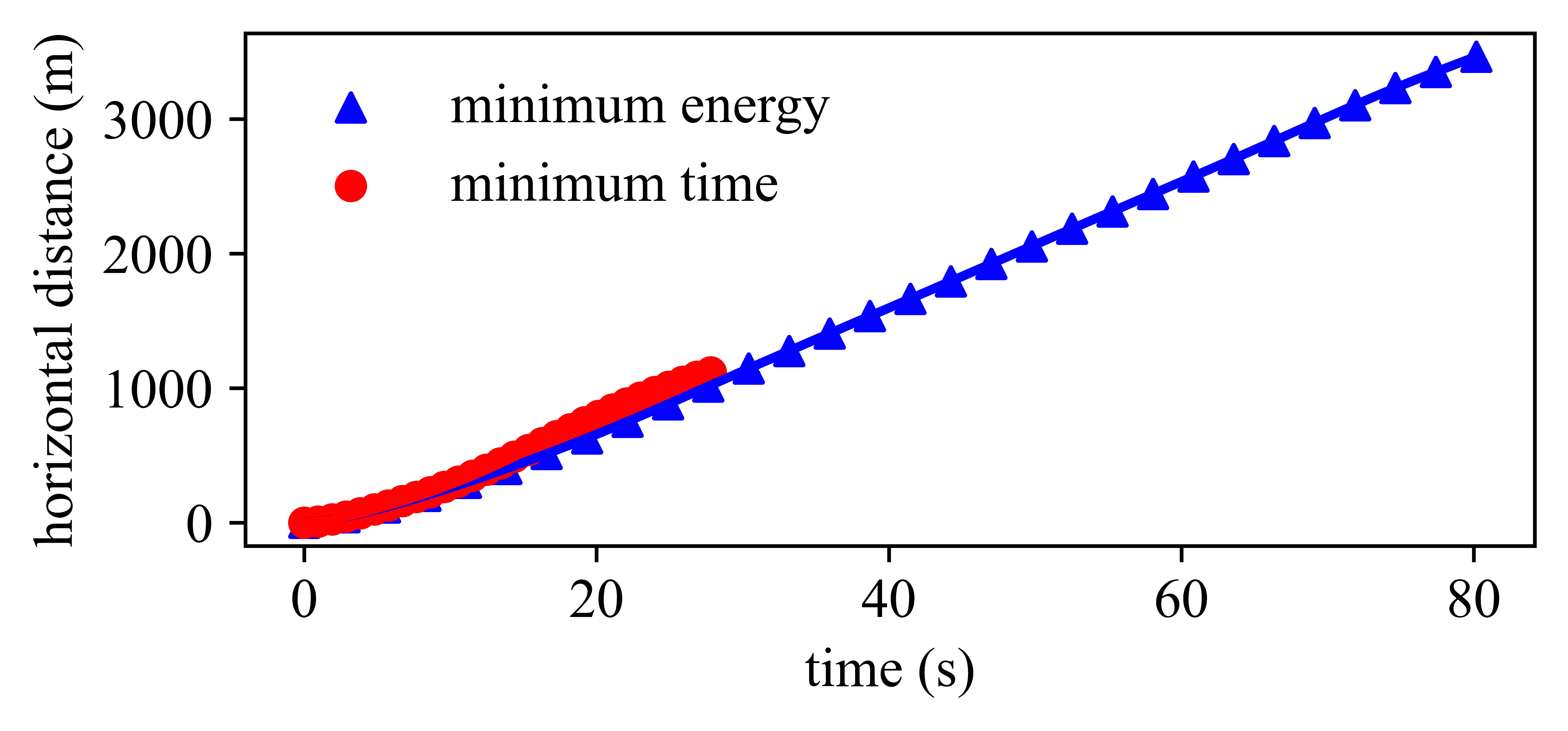}
\end{subfigure}
\caption{Time-histories of aircraft states $\mathbf{(v,\gamma,h,x)}$ for minimum energy and minimum time objectives.}
\end{figure}

The aircraft accelerates past the target velocity for both objectives (Fig. 8). The velocity begins to decrease once the aircraft begins to climb, and only returns to the target velocity at the termination of the trajectory. This makes sense for minimum time trajectories, which aim to go as fast as possible to minimize total time. This is less obvious for the minimum energy case. For these trajectories we see that the aircraft climbs at roughly a $1^{\circ}$ angle of attack in order to minimize drag (see Fig. 4). Because of this, the aircraft must go faster than expected in order to generate enough lift to climb. The minimum energy trajectories also exhibit an increased rate of climb near the end of the transition. This bleeds off the excess velocity, thereby enabling the aircraft to reach the final velocity constraint at the same time as the final altitude constraint.

The profiles of flight path angle vary significantly between objectives. Minimum energy trajectories maintain a flight path angle near zero while the aircraft accelerates horizontally, before the aircraft pitches upwards to begin climbing to the target altitude. The minimum time trajectory contains an oscillatory flight path angle profile, which roughly coincides with the crossover point identified earlier between lift generated by the rotors and lift generated by the wing. The flight path angle remains positive at the end of the minimum time trajectory, indicating the aircraft has a positive rate of climb at the target altitude. This contrasts with the flight path angle for the minimum energy objective which returns to near-zero towards the end of the transition, indicating that the aircraft levels off at the target altitude.

\FloatBarrier
\begin{figure}[!h]
\centering
\begin{subfigure}[h]{0.45\linewidth}
\includegraphics[width=1\linewidth]{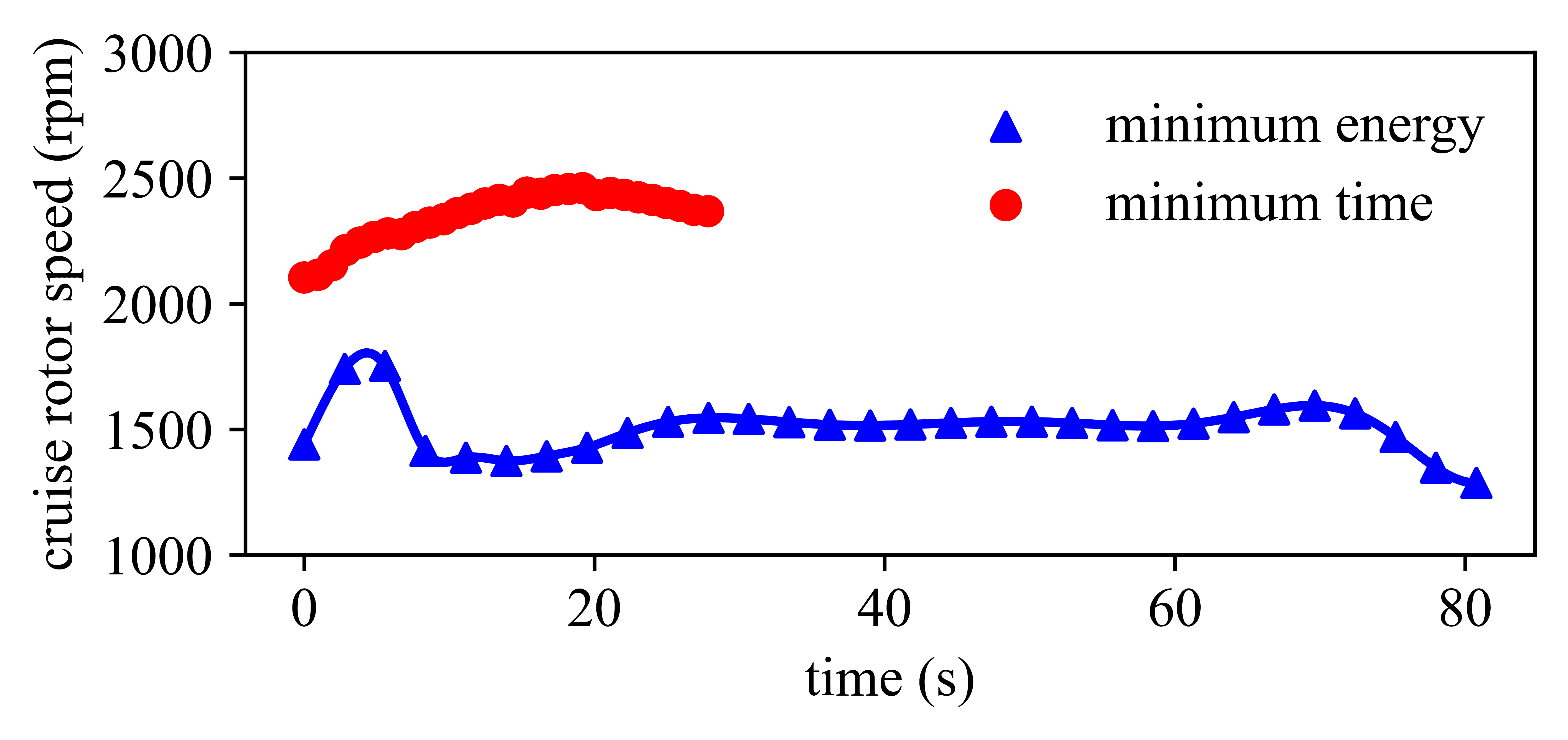}
\end{subfigure}
\begin{subfigure}[h]{0.45\linewidth}
\includegraphics[width=1\linewidth]{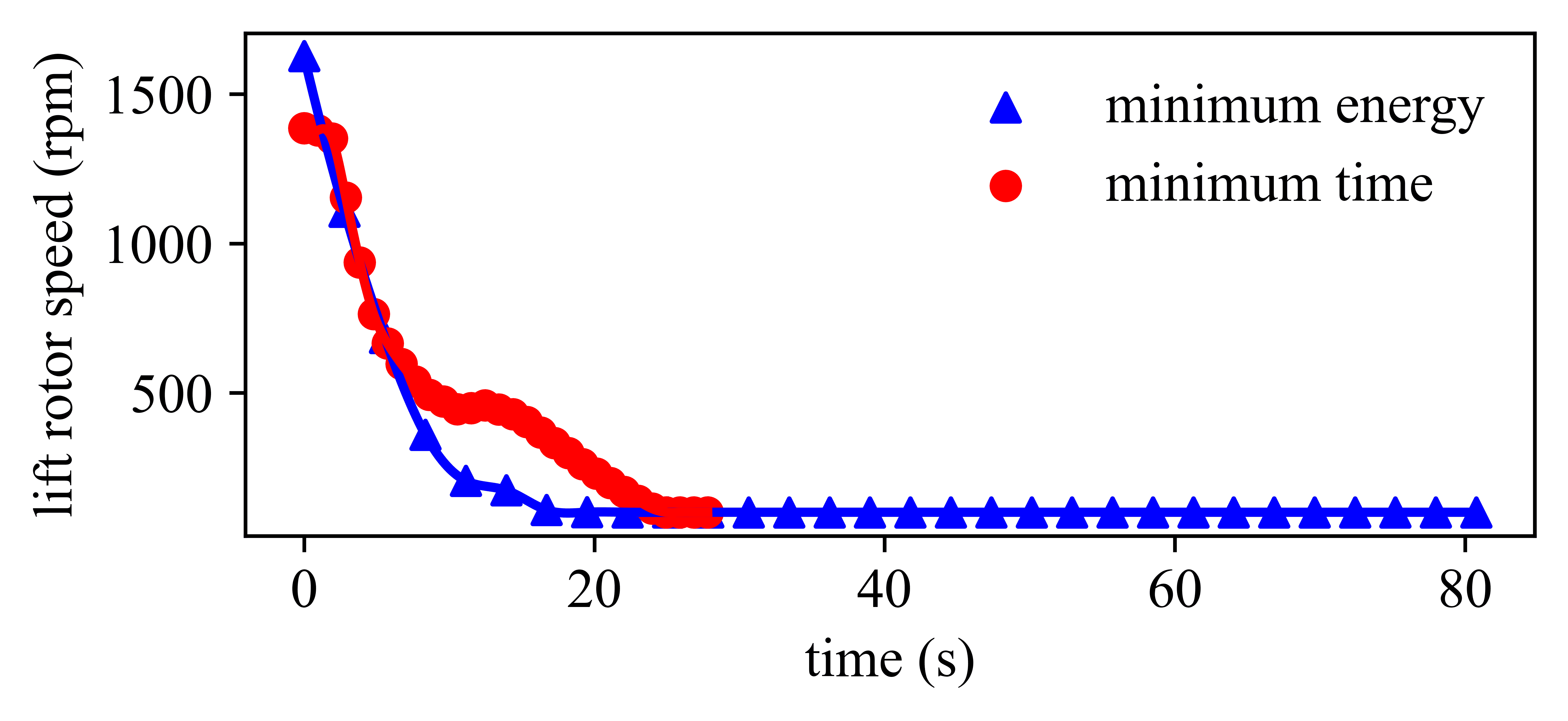}
\end{subfigure}
\begin{subfigure}[h]{0.45\linewidth}
\includegraphics[width=1\linewidth]{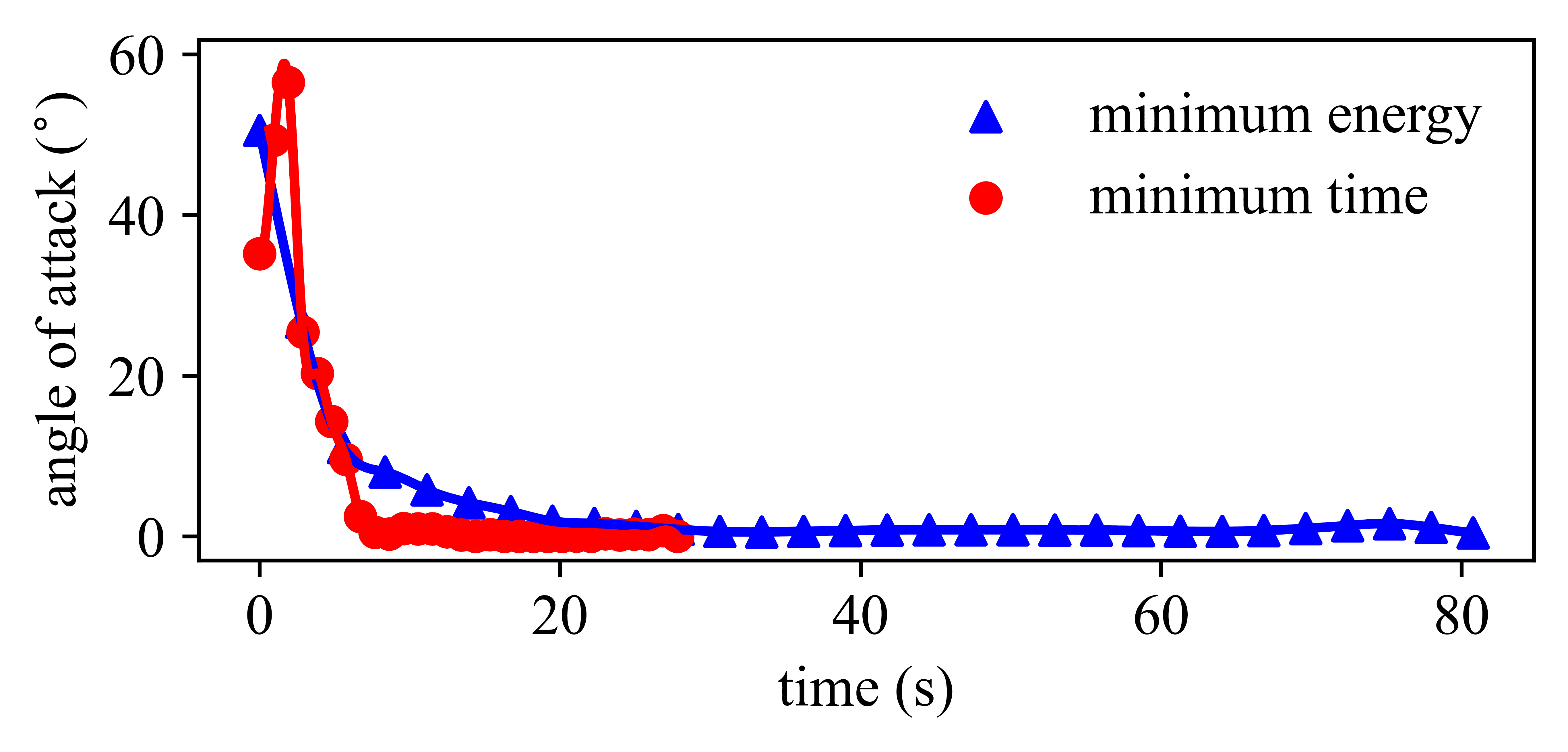}
\end{subfigure}
\begin{subfigure}[h]{0.45\linewidth}
\includegraphics[width=1\linewidth]{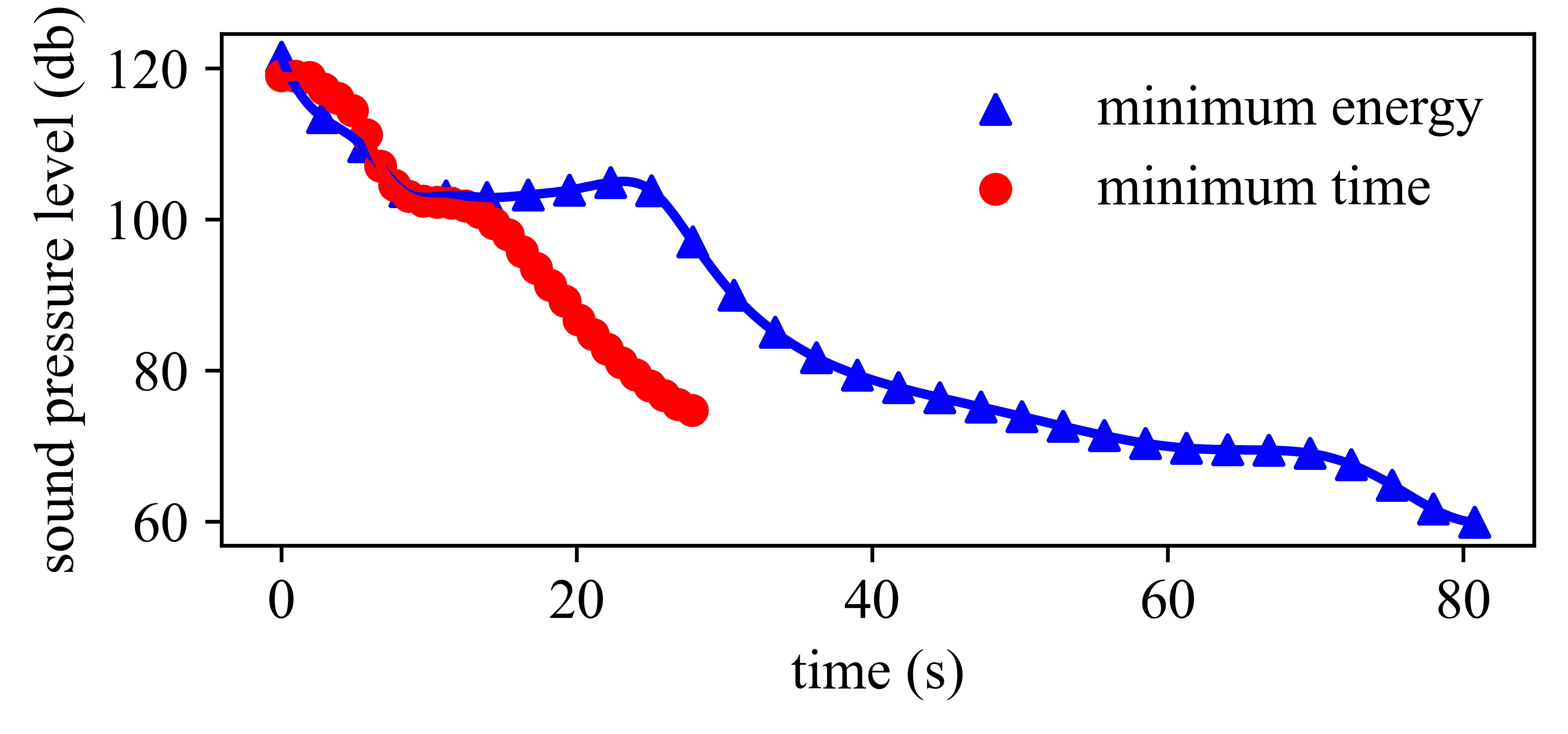}
\end{subfigure}
\caption{Cruise rotor speed, lift rotor speed, angle of attack, and sound pressure level for minimum energy and minimum time objectives.}
\end{figure}

\begin{figure}[ht]
\centering
\includegraphics[width=0.9\textwidth]{e_vs_t.png}
\caption{Energy expended versus time for minimum energy and minimum time trajectories.}
\end{figure}

\begin{table}[ht]
\begin{spacing}{0.85}
\caption{\label{nom} Results for minimum energy and minimum time trajectories}
\centering
\begin{tabular}{ r c c c }
\hline \hline
& Total energy (MJ / Wh) & Transition time (s) & Transition efficiency (\%) \\ 
Minimum energy (nominal) & 13.29 / 3691 & 80.16 & 33.0\\
Minimum time (nominal) & 16.35 / 4542 & 27.81 & 26.8\\
\hline
Minimum energy $(|\theta| \leq 45^{\circ})$ & 13.32 / 3692 & 83.19 & 32.9\\
Minimum energy $(|\theta| \leq 30^{\circ})$ & 13.54 / 3761 & 86.27 & 32.4\\
Minimum energy $(|\theta| \leq 15^{\circ})$ & 13.67 / 3797 & 122.94 & 32.1\\
\hline
Minimum time $(|\theta| \leq 45^{\circ})$ & 16.10 / 4472 & 27.95 & 27.2\\
Minimum time $(|\theta| \leq 30^{\circ})$ & 16.94 / 4706 & 29.85 & 27.4\\
Minimum time $(|\theta| \leq 15^{\circ})$ & 26.28 / 7300 & 46.86 & 27.3\\
\hline \hline
\end{tabular}
\end{spacing}
\end{table}

\FloatBarrier
\subsection{Pitch Angle Constraints}

The maximum pitch angles ($\theta=\gamma+\alpha$) for the minimum energy and minimum time trajectories are $51^{\circ}$ and $56^{\circ}$ respectively. These angles are likely to be unacceptable for a UAM passenger service, where passenger comfort is a priority (the typical maximum airliner pitch angle is around 15-20$^{\circ}$~\cite{wakefield2009exceeding}). Figures 11 and 12 show the resulting trajectories with maximum pitch angle constraints added to the optimization problem.

Qualitatively, very little difference can be seen between the transition trajectories for both objectives with $45^{\circ}$ and $30^{\circ}$ pitch angle constraints. Quantitatively, these trajectories perform only slightly worse than their unconstrained counterparts (Table 2) with both objectives seeing a change in efficiency of less than $1\%$ for $45^{\circ}$ and $30^{\circ}$ constraints. This small efficiency difference is representative of a very flat energy design-space with respect to maximum pitch angle. This is similar to the conclusions of Chauhan and Martins~\cite{chauhan2020tilt}. Additionally, the transition time increases by nearly 8\% for minimum energy trajectories with $30^{\circ}$ constraints, and less than 1\% for minimum time trajectories.

\FloatBarrier
\begin{figure}[ht]
\centering
\includegraphics[width=0.9\textwidth]{new_e_theta_h_vs_x.png}
\caption{Altitude versus horizontal displacement for minimum energy trajectories with pitch angle constraints.}
\end{figure}

\FloatBarrier
\begin{figure}[ht]
\centering
\includegraphics[width=0.9\textwidth]{new_t_theta_h_vs_x.png}
\caption{Altitude versus horizontal displacement for minimum time trajectories with pitch angle constraints.}
\end{figure}

The trajectories constrained by $|\theta|\leq15^{\circ}$ show substantial differences. For both objectives these trajectories cover more horizontal distance (5406m and 2458m for minimum energy and minimum time respectively). This represents a 146\% and 204\% increase when compared to the $30^{\circ}$ constraints. Of particular interest is the increasing similarity between the results for both objectives as the pitch constraints become more aggressive. Figure 13 shows that the pitch angle profiles exhibit similarity for the $15^{\circ}$ constraint. Despite this visual similarity, the energy consumption is drastically different, with minimum time trajectories using nearly two times the energy of the minimum energy trajectories for the $15^{\circ}$ constraint.

\FloatBarrier
\begin{figure}[ht]
\centering
\begin{subfigure}[h]{0.45\linewidth}
\includegraphics[width=1\linewidth]{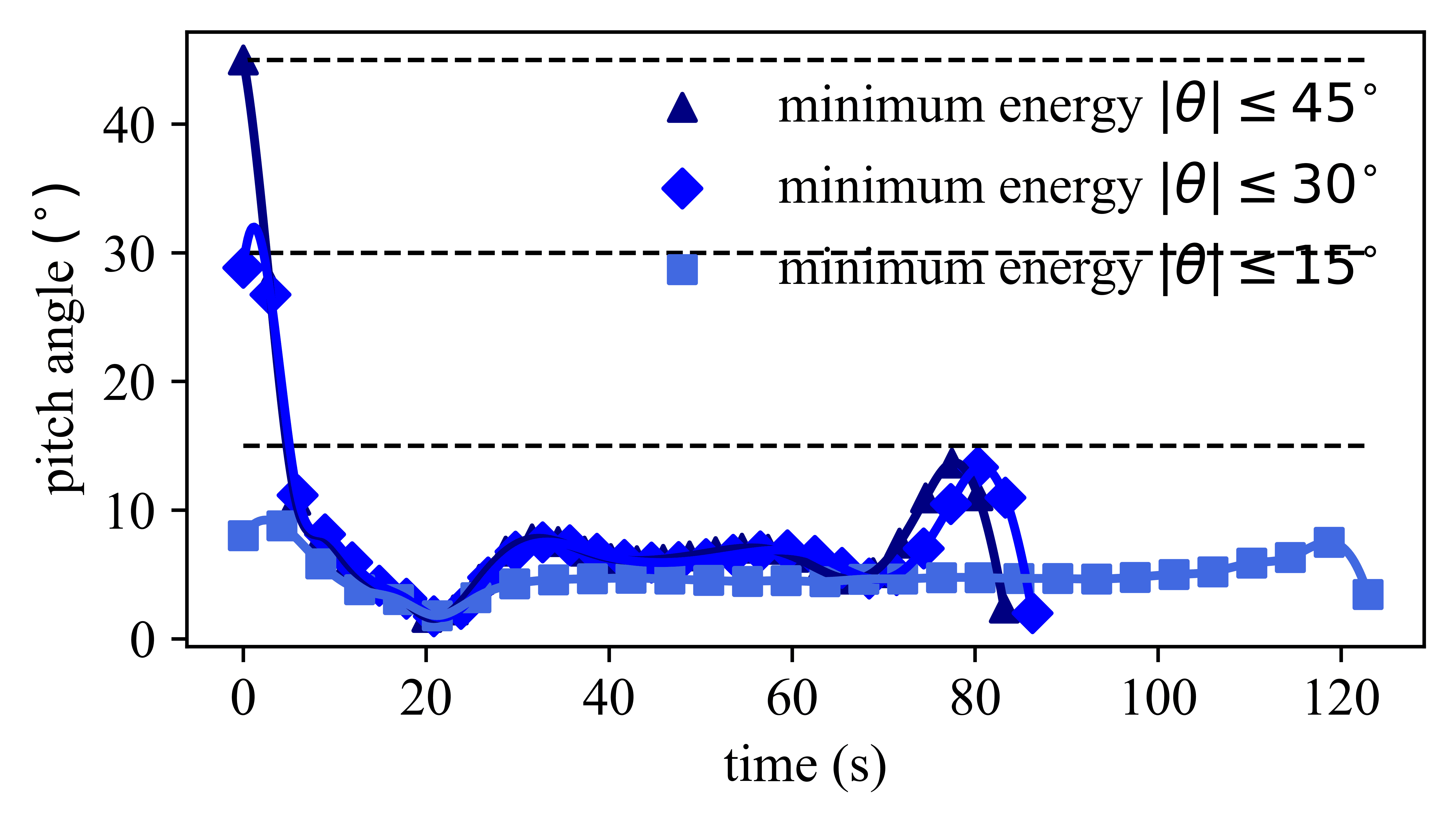}
\end{subfigure}
\begin{subfigure}[h]{0.45\linewidth}
\includegraphics[width=1\linewidth]{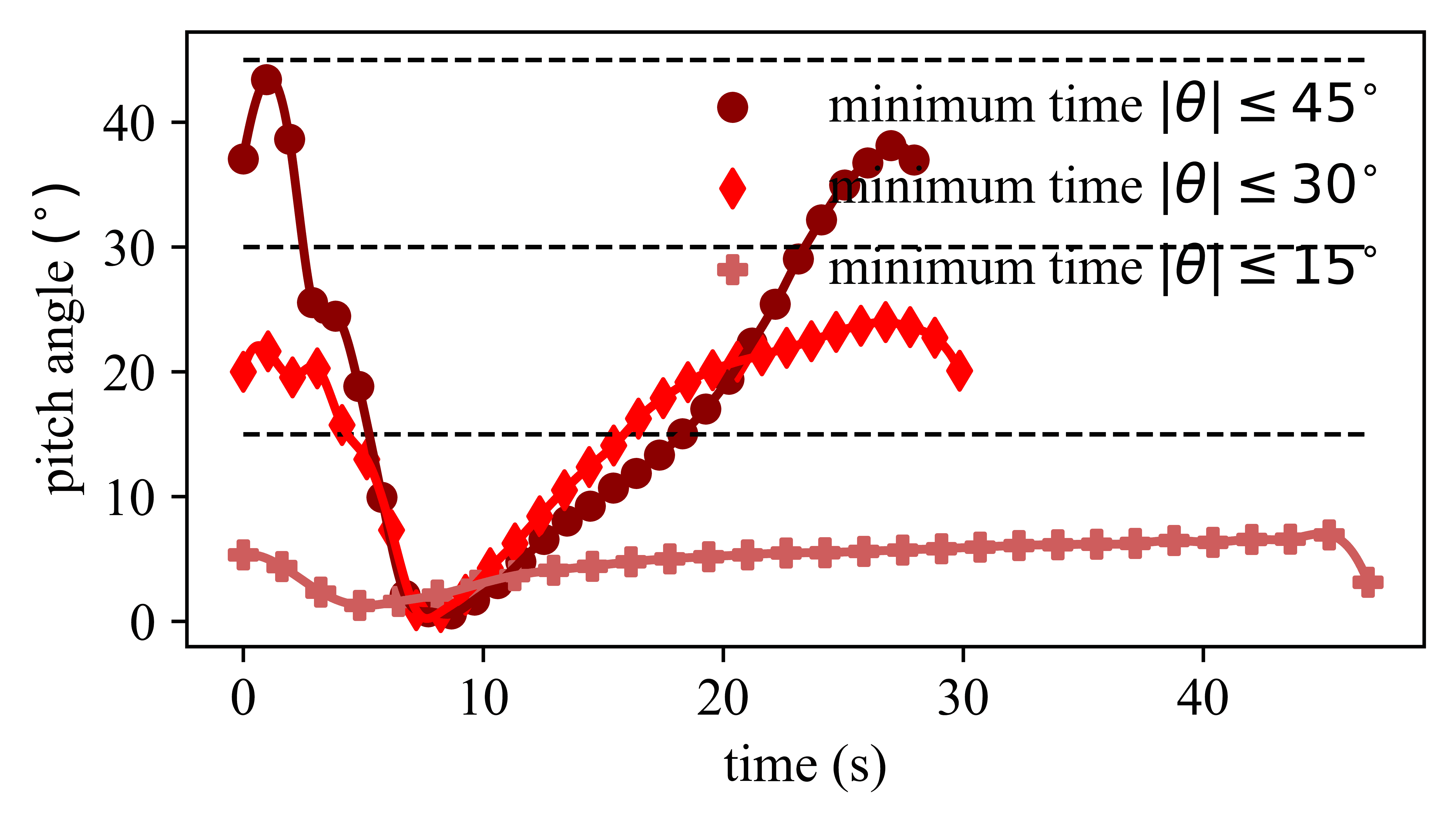}
\end{subfigure}
\caption{Pitch angle versus time for $\mathbf{45^{\circ}}$, $\mathbf{30^{\circ}}$, and $\mathbf{15^{\circ}}$ constraints.}
\end{figure}

\subsection{Acoustic Constraints}

Implementing acoustic constraints presents certain practical difficulties. At the start of the transition maneuver, aircraft velocity is low and the lifting rotors are the only source of lift. Since the minimum altitude is constrained $(\min{\vec{h}}\geq h_0)$, these rotors are forced to provide sufficient thrust to sustain the aircraft altitude. This means that if we constrain the initial rotor noise we risk creating an infeasible optimization problem since the aircraft cannot satisfy the minimum altitude constraint and the maximum noise constraint on takeoff. Effectively, the initial rotor noise has a minimum value, beyond which the problem becomes infeasible. 

We solve this problem by evaluating two different types of acoustic constraints. First, we apply a constant noise constraint to the final 2/3 (as a function of time) of the transition trajectory. This allows the aircraft to takeoff without noise constraints, thereby removing constraints on the loudest portion of the transition. Second, we apply a linear acoustic constraint. We constrain the initial SPL by a large value (120db) and linearly reduce the constraint from $t_0$ to $t_f$. This method simulates the effect of reducing overall sound pressure across the entire transition maneuver while still allowing for a noisy takeoff. Both these methods provide valuable insight into the characteristics of noise limited trajectories while maintaining a feasible optimization problem.

\FloatBarrier
\begin{figure}[ht]
\centering
\includegraphics[width=0.9\textwidth]{new_e_spl_h_vs_x.png}
\caption{Altitude versus horizontal displacement for minimum energy trajectories with noise constraints over the final 2/3 of the trajectory.}
\end{figure}

\FloatBarrier
\begin{figure}[ht]
\centering
\includegraphics[width=0.9\textwidth]{new_t_spl_h_vs_x.png}
\caption{Altitude versus horizontal displacement for minimum time trajectories with noise constraints over the final 2/3 of the trajectory.}
\end{figure}

Figures 14 and 15 show  the results of the first acoustic constraint for both objectives. In these plots, the ground-level SPL is constrained by a constant over the latter 2/3 of the trajectory. Clearly, minimum energy trajectories are only slightly affected by these constraints. In fact the only visual difference between these trajectories is at the inflection point where the aircraft begins to climb. At this point the aircraft gains more initial altitude in order to reduce the ground-level SPL.

The minimum time trajectories are significantly more affected. By constraining the latter 2/3 of the trajectory we force the trajectories towards aggressive altitude and velocity gains in the first 1/3 of the trajectory. This is because rotor noise is correlated with rotor power, and increased rotor power tends to reduce transition time, thereby resulting in optimal minimum time trajectories.

This is very apparent for the $SPL\leq75$db constraint where the aircraft immediately climbs to an altitude near the target altitude and then slowly accelerates to the target velocity over a large amount of time. The rotor noise corresponding with these trajectories is shown in Fig. 16. We remind the reader that the decibel scale is nonlinear, and that linear differences between trajectories are not expected (despite the linearly decreasing nature of the constraints).

\FloatBarrier
\begin{figure}[ht]
\centering
\begin{subfigure}[h]{0.45\linewidth}
\includegraphics[width=1\linewidth]{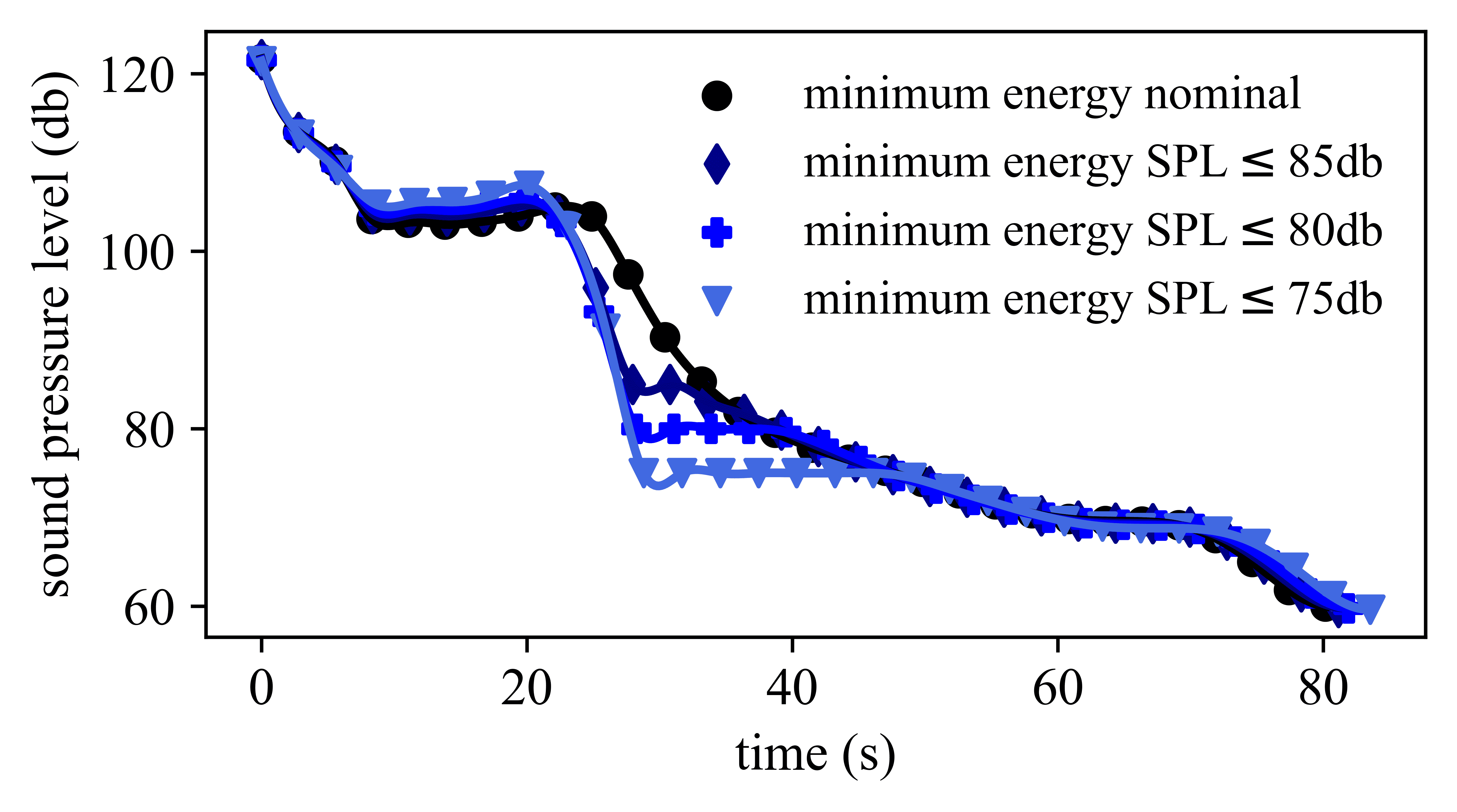}
\end{subfigure}
\begin{subfigure}[h]{0.45\linewidth}
\includegraphics[width=1\linewidth]{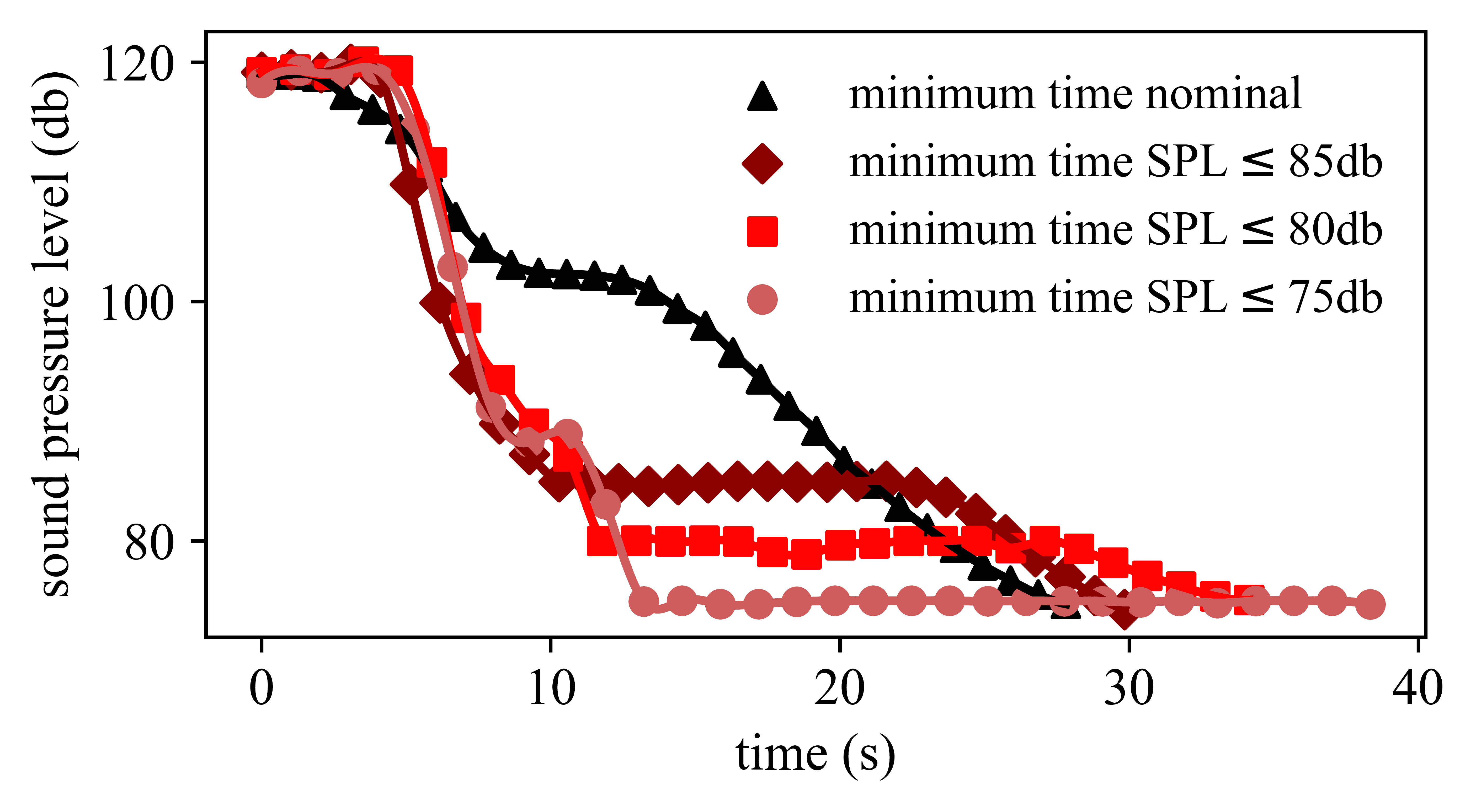}
\end{subfigure}
\caption{Sound pressure level versus time for 85db, 80db, and 75db constraints over the final 2/3 of the trajectory.}
\end{figure}

Table 4 summarizes the results for the noise-constrained trajectories. We find that transition efficiency is significantly decreased for minimum time trajectories (approximately 0.33\% per db) and transition time is (approximately 0.6s per db). From this we conclude that acoustic constraints have a large effect on minimum time trajectories and that the transition time design space is not flat with respect to rotor noise (unlike the results with pitch constraints). We also see that the energy cost of noise-constrained minimum time trajectories is quite high, with an increase of nearly 6MJ across all minimum time trajectories.

The noise constraints used here are all greater than the 67db at 250ft determined by Uber Elevate~\cite{holden2016fast} as a requirement for urban air mobility. In fact, the Lift-plus-Cruise air taxi (as modeled in this paper) cannot perform a complete transition under this constraint. This presents a strong case for large-scale simultaneous optimization of both trajectory and aircraft design.

\begin{table}[ht]
\begin{spacing}{0.85}
\caption{\label{acoustic} Results for minimum energy and minimum time trajectories with noise constraints}
\centering
\begin{tabular}{ r c c c }
\hline \hline
& Total energy (MJ / Wh) & Transition time (s) & Transition efficiency (\%) \\ 
\hline
Minimum energy $(SPL_{2/3} \leq \text{85db})$ & 13.30 / 3694 & 81.13 & 32.9\\
Minimum energy $(SPL_{2/3} \leq \text{80db})$ & 13.32 / 3700 & 81.89 & 32.9\\ 
Minimum energy $(SPL_{2/3} \leq \text{75db})$ & 13.36 / 3711 & 83.53 & 32.8\\
\hline
Minimum time $(SPL_{2/3} \leq \text{85db})$ & 19.16 / 5322 & 29.84 & 22.9\\
Minimum time $(SPL_{2/3} \leq \text{80db})$ & 21.79 / 6053 & 34.15 & 20.1\\
Minimum time $(SPL_{2/3} \leq \text{75db})$ & 24.79 / 6886 & 38.34 & 18.0\\
\hline
Minimum energy $(SPL_{lin} \leq \text{120-55db})$ & 13.36 / 3711 & 83.39 & 32.8\\
Minimum energy $(SPL_{lin} \leq \text{120-50db})$ & 13.44 / 3733 & 84.85 & 32.7\\
Minimum energy $(SPL_{lin} \leq \text{120-45db})$ & 13.56 / 3767 & 86.56 & 32.4\\
\hline
Minimum time $(SPL_{lin} \leq \text{120-70db})$ & 16.66 / 4628 & 30.41 & 26.3\\
Minimum time $(SPL_{lin} \leq \text{120-65db})$ & 17.15 / 4764 & 33.09 & 25.6\\
Minimum time $(SPL_{lin} \leq \text{120-60db})$ & 17.35 / 4819 & 35.43 & 25.3\\
\hline \hline
\end{tabular}
\end{spacing}
\end{table}

Figures 17-19 depict the results of the linearly varying noise constraint. We examine three different linear interpolations for minimum energy trajectories: 120-55db, 120-50db, and 120-45db. For minimum time trajectories we us different values: 120-70db, 120-65db, and 120-60db. This is because the minimum energy trajectories would be altogether unaffected by the higher constraints used for minimum time, and minimum time trajectories would be significantly affected by any lower constraints. This would make comparisons between nominal trajectories and constrained trajectories difficult.



\FloatBarrier
\begin{figure}[!ht]
\centering
\includegraphics[width=0.9\textwidth]{new_e_lin_spl_h_vs_x.png}
\caption{Altitude versus horizontal displacement for minimum energy trajectories with linear noise constraints.}
\end{figure}

\FloatBarrier
\begin{figure}[!ht]
\centering
\includegraphics[width=0.9\textwidth]{new_t_lin_spl_h_vs_x.png}
\caption{Altitude versus horizontal displacement for minimum time trajectories with linear noise constraints.}
\end{figure}

\FloatBarrier
Minimum energy trajectories are affected very little by all the linear constraints with an efficiency difference of only 0.4\% across all linear constraints and virtually no difference in transition time (2.4s). This is explained by Figure 19, which shows that the acoustic constraints are only active for a small portion of each trajectory (at the beginning of the climb and at the end of the transition). This is because the SPL for the nominal trajectories already follows a roughly linear descent, and further constraints have minimal effect. Minimum time trajectories are more affected. There is a 1\% difference in transition efficiency across all the linear constraints and the trajectories noticeably exhibit different behavior. Specifically, the 120-60db trajectory climbs to the target altitude before accelerating horizontally (similar to the behavior with constant noise constraints on the latter 2/3 of the transition).

Figure 18 shows an interesting result for minimum time trajectories. The slope of all constraints and rotor noise curves are roughly identical. This implies that optimal, minimum time trajectories have a strong correlation with maximum noise trajectories, and that there is very little difference in the noise profiles of these trajectories.

\FloatBarrier
\begin{figure}[ht]
\centering
\begin{subfigure}[h]{0.45\linewidth}
\includegraphics[width=1\linewidth]{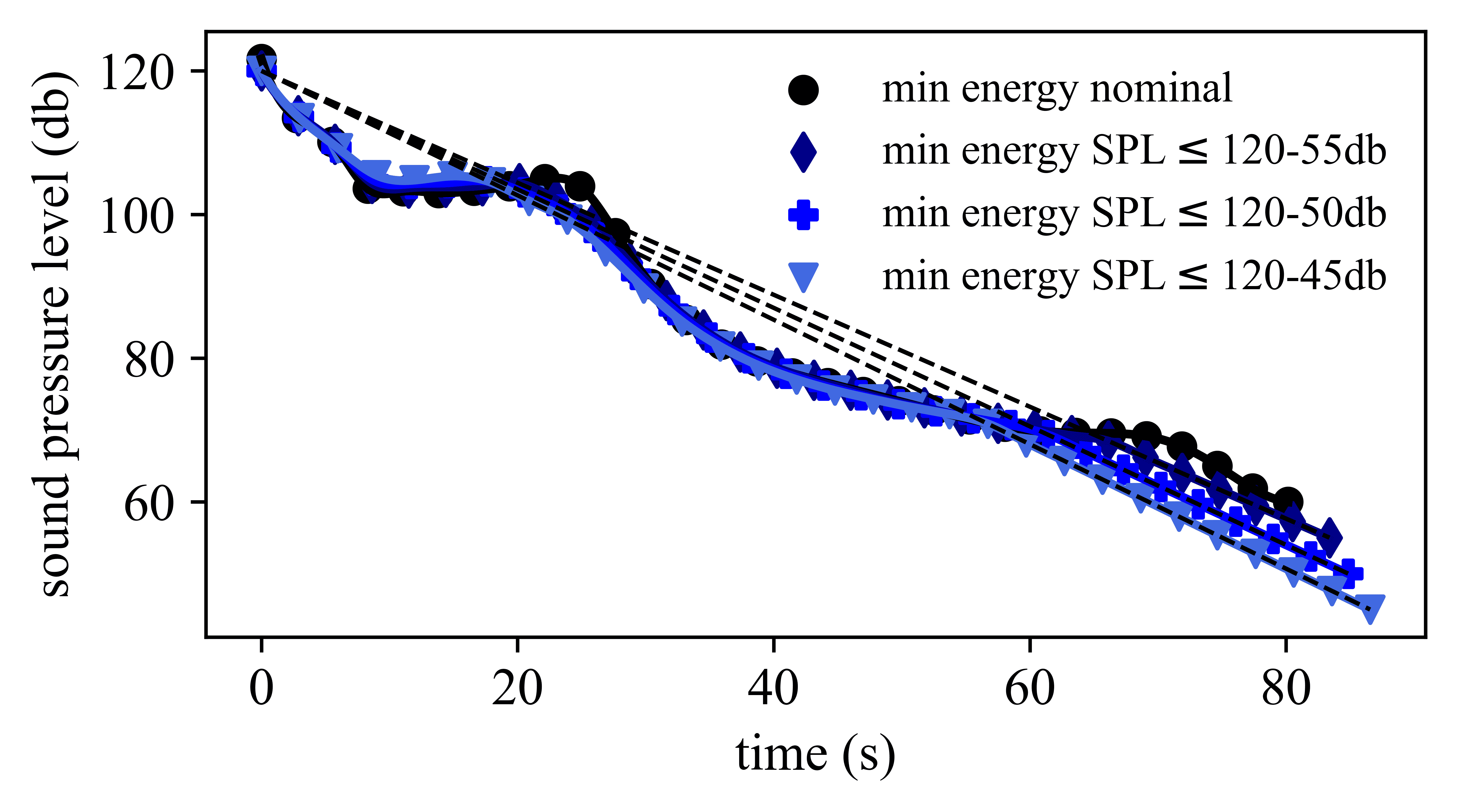}
\end{subfigure}
\begin{subfigure}[h]{0.45\linewidth}
\includegraphics[width=1\linewidth]{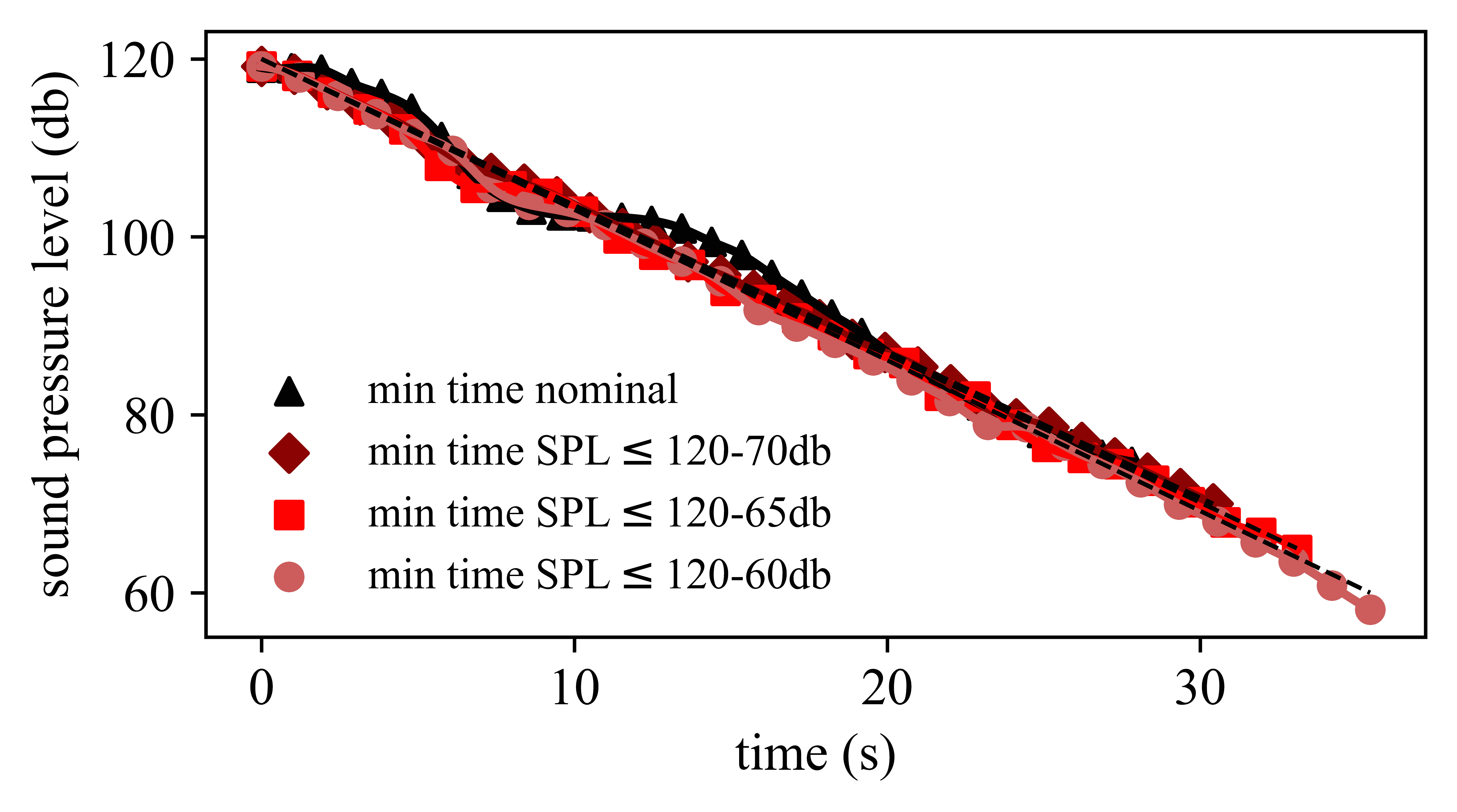}
\end{subfigure}
\caption{Sound pressure level versus time for linear noise constraints.}
\end{figure}

\section{Conclusion}

In this paper we present a trajectory optimization method that combines models for flight dynamics, rotor aerodynamics, wing aerodynamics, and motor performance, and makes use of direct transcription and large-scale MDO. A principle contribution of this paper is this combination of methods which proved to be highly effective.

We demonstrate the versatility of this combination of methods by finding optimal transition trajectories for NASA's Lift-plus-Cruise air taxi that satisfy minimum energy and minimum time objectives. Results for these objectives are compared, along with trajectories that are generated using constraints on pitch angle and noise.

We find several notable differences between unconstrained minimum energy and minimum time trajectories. First, we note that minimum energy trajectories are characterized by two segments: constant-altitude acceleration followed by a second segment with a steady rate of climb. This contrasts with minimum-time trajectories which immediately begin to climb and exhibit a noticeable trajectory change when switching from lift generated by the rotors to lift generated by the wing. 

Second, we find that minimum-time trajectories travel less than half the horizontal distance of minimum energy trajectories. As expected, this is due to increased cruise rotor power which results in a significantly higher rate of climb. Despite the visual differences in these trajectories, the transition efficiency is only 6.2\% different (without additional constraints) which is representative of a very flat energy space with respect to both objectives. Interestingly, the lift rotors play a very small role during transition for both objectives. In fact, lift rotor speed and power decrease rapidly to near zero, and remain low throughout the duration of the transition.

We evaluate the effect of pitch angle constraints on transition trajectories for both objectives. We find that trajectories are relatively unchanged for $|\theta|\leq45^{\circ}$ and $|\theta|\leq30^{\circ}$ pitch angle constraints. However, the $|\theta|\leq15^{\circ}$ constraints result in larger transition times and increased horizontal distance traveled. In general, we show that $|\theta|\leq15^{\circ}$ constraints reduce the rate of climb, resulting in longer transition times.

As expected, aircraft noise is a barrier to UAM adoption. Therefore we also evaluate two types of acoustic constraints. First, we constrain the latter two thirds of the trajectory by constant sound pressure limits of 85db, 80db, and 75db. We find that the minimum energy trajectories are relatively unaffected by this constraint, whereas minimum time trajectories can be generally characterized by a period of more rapid initial altitude gain followed by a period of slower horizontal acceleration. Second, we evaluate a linearly reducing acoustic constraint across the duration of the transition and once again find that minimum energy trajectories are less affected than minimum time trajectories. However, both sets of trajectories change very little compared to unconstrained results (0.6\% and 1.5\% maximum efficiency differences for minimum energy and minimum time trajectories respectively). 

Across all acoustic constraints we find, as expected, that the aircraft prioritizes altitude gain and acceleration in the less-constrained stages of the transition. The models used here and our analysis suggests that the Lift-plus-Cruise air taxi, in its current configuration, might exceed the proposed limits for UAM noise levels during the transition trajectory. This presents an opportunity to use large-scale MDO to simultaneously optimize the aircraft design and the transition trajectory. (The current study only considers broadband noise and not, for instance, tonal noise which might constrain the trajectories even more).

Overall, we demonstrate the application of direct transcription methods to transition trajectory optimization for the Lift-plus-Cruise aircraft. In doing so, we show the versatility of this approach by including several physics-based models and constraints. Our results are practical, and reinforce the need for the analysis of dynamic mission segments in the context of aircraft design by demonstrating, among other things, the changes that pitch angle and acoustic constraints can have on trajectories.

\section*{Acknowledgments}
The material presented in this paper is, in part, based upon work supported by NASA under award No.~80NSSC21M0070 and by DARPA under grant No. D23AP00028-00.

\bibliography{main}

\begin{thebibliography}{49}
\newcommand{\enquote}[1]{``#1''}
\providecommand{\natexlab}[1]{#1}
\providecommand{\url}[1]{\texttt{#1}}
\providecommand{\urlprefix}{URL }
\expandafter\ifx\csname urlstyle\endcsname\relax
  \providecommand{\doi}[1]{\discretionary{}{}{}https://doi.org/#1}\else
  \providecommand{\doi}[1]{\discretionary{}{}{}\urlstyle{rm}\url{https://doi.org/#1}}\fi

\bibitem[{Hwang and Ning(2018)}]{hwang2018large}
Hwang, J.~T., and Ning, A., \enquote{Large-scale multidisciplinary optimization
  of an electric aircraft for on-demand mobility,} \emph{2018 AIAA/ASCE/AHS/ASC
  Structures, Structural Dynamics, and Materials Conference}, 2018, p. 1384.
\newblock \doi{10.2514/6.2018-1384}.

\bibitem[{Silva et~al.(2018)Silva, Johnson, Solis, Patterson, and
  Antcliff}]{silva2018vtol}
Silva, C., Johnson, W.~R., Solis, E., Patterson, M.~D., and Antcliff, K.~R.,
  \enquote{VTOL urban air mobility concept vehicles for technology
  development,} \emph{2018 Aviation Technology, Integration, and Operations
  Conference}, 2018, p. 3847.

\bibitem[{Chauhan and Martins(2020)}]{chauhan2020tilt}
Chauhan, S.~S., and Martins, J.~R., \enquote{Tilt-wing eVTOL takeoff trajectory
  optimization,} \emph{Journal of aircraft}, Vol.~57, No.~1, 2020, pp. 93--112.
\newblock \doi{10.2514/1.C035476}.

\bibitem[{Kubo and Suzuki(2008)}]{kubo2008tail}
Kubo, D., and Suzuki, S., \enquote{Tail-sitter vertical takeoff and landing
  unmanned aerial vehicle: transitional flight analysis,} \emph{Journal of
  Aircraft}, Vol.~45, No.~1, 2008, pp. 292--297.
\newblock \doi{10.2514/1.30122}.

\bibitem[{Anderson et~al.(2021)Anderson, Willis, Johnson, Ning, and
  Beard}]{anderson2021comparison}
Anderson, R., Willis, J., Johnson, J., Ning, A., and Beard, R.~W., \enquote{A
  comparison of aerodynamics models for optimizing the takeoff and transition
  of a bi-wing tailsitter,} \emph{AIAA Scitech 2021 Forum}, 2021, p. 1008.
\newblock \doi{10.2514/6.2021-1008}.

\bibitem[{Betts(1998)}]{betts1998survey}
Betts, J.~T., \enquote{Survey of numerical methods for trajectory
  optimization,} \emph{Journal of guidance, control, and dynamics}, Vol.~21,
  No.~2, 1998, pp. 193--207.
\newblock \doi{10.2514/2.4231}.

\bibitem[{Kelly(2017)}]{kelly2017}
Kelly, M., \enquote{An introduction to trajectory optimization: How to do your
  own direct collocation,} \emph{SIAM Review}, Vol.~59, No.~4, 2017, pp.
  849--904.
\newblock \doi{10.1137/16M1062569}.

\bibitem[{Sutton and Barto(2018)}]{sutton2018reinforcement}
Sutton, R.~S., and Barto, A.~G., \emph{Reinforcement learning: An
  introduction}, MIT press, 2018.

\bibitem[{Hargraves and Paris(1987)}]{hargraves1987direct}
Hargraves, C.~R., and Paris, S.~W., \enquote{Direct trajectory optimization
  using nonlinear programming and collocation,} \emph{Journal of guidance,
  control, and dynamics}, Vol.~10, No.~4, 1987, pp. 338--342.
\newblock \doi{10.2514/3.20223}.

\bibitem[{Bryson and Denham(1962)}]{bryson1962steepest}
Bryson, A.~E., and Denham, W.~F., \enquote{A steepest-ascent method for solving
  optimum programming problems,} \emph{Journal of Applied Mechanics}, 1962.

\bibitem[{Hwang and Munster(2018)}]{hwang2018solution}
Hwang, J.~T., and Munster, D., \enquote{Solution of ordinary differential
  equations in gradient-based multidisciplinary design optimization,}
  \emph{2018 AIAA/ASCE/AHS/ASC Structures, Structural Dynamics, and Materials
  Conference}, 2018, p. 1646.

\bibitem[{Martins and Hwang(2013)}]{martins2013review}
Martins, J.~R., and Hwang, J.~T., \enquote{Review and unification of methods
  for computing derivatives of multidisciplinary computational models,}
  \emph{AIAA journal}, Vol.~51, No.~11, 2013, pp. 2582--2599.
\newblock \doi{10.2514/1.J052184}.

\bibitem[{Hwang and Martins(2018{\natexlab{a}})}]{hwang2018computational}
Hwang, J.~T., and Martins, J.~R., \enquote{A computational architecture for
  coupling heterogeneous numerical models and computing coupled derivatives,}
  \emph{ACM Transactions on Mathematical Software (TOMS)}, Vol.~44, No.~4,
  2018{\natexlab{a}}, pp. 1--39.
\newblock \doi{10.1145/3182393}.

\bibitem[{Betts and Cramer(1995)}]{betts1995application}
Betts, J.~T., and Cramer, E.~J., \enquote{Application of direct transcription
  to commercial aircraft trajectory optimization,} \emph{Journal of Guidance,
  Control, and Dynamics}, Vol.~18, No.~1, 1995, pp. 151--159.
\newblock \doi{10.2514/3.56670}.

\bibitem[{Hargraves et~al.(1981)Hargraves, Johnson, Paris, and
  Retties}]{hargraves1981numerical}
Hargraves, C., Johnson, F., Paris, S., and Retties, I., \enquote{Numerical
  computation of optimal atmospheric trajectories,} \emph{Journal of Guidance
  and Control}, Vol.~4, No.~4, 1981, pp. 406--414.
\newblock \doi{10.2514/3.56093}.

\bibitem[{Hendricks et~al.(2019)Hendricks, Falck, Gray, Aretskin-Hariton,
  Ingraham, Chapman, Schnulo, Chin, Jasa, and
  Bergeson}]{hendricks2019multidisciplinary}
Hendricks, E.~S., Falck, R.~D., Gray, J.~S., Aretskin-Hariton, E., Ingraham,
  D., Chapman, J.~W., Schnulo, S.~L., Chin, J., Jasa, J.~P., and Bergeson,
  J.~D., \enquote{Multidisciplinary optimization of a turboelectric tiltwing
  urban air mobility aircraft,} \emph{AIAA Aviation 2019 Forum}, 2019, p. 3551.
\newblock \doi{10.2514/6.2019-3551}.

\bibitem[{Jasa et~al.(2020)Jasa, Brelje, Gray, Mader, and
  Martins}]{jasa2020large}
Jasa, J.~P., Brelje, B.~J., Gray, J.~S., Mader, C.~A., and Martins, J.~R.,
  \enquote{Large-Scale Path-Dependent Optimization of Supersonic Aircraft,}
  \emph{Aerospace}, Vol.~7, No.~10, 2020, p. 152.
\newblock \doi{10.3390/aerospace7100152}.

\bibitem[{Falck et~al.(2018)Falck, Ingraham, and
  Aretskin-Hariton}]{falck2018multidisciplinary}
Falck, R.~D., Ingraham, D., and Aretskin-Hariton, E.,
  \enquote{Multidisciplinary optimization of urban-air-mobility class Aircraft
  trajectories with acoustic constraints,} \emph{2018 AIAA/IEEE Electric
  Aircraft Technologies Symposium (EATS)}, IEEE, 2018, pp. 1--7.
\newblock \doi{10.2514/6.2018-4985}.

\bibitem[{Orndorff and Hwang(2022)}]{orndorff2022investigation}
Orndorff, N.~C., and Hwang, J.~T., \enquote{Investigation of Optimal Air-Taxi
  Transition Profiles using Direct-Transcription Trajectory Optimization,}
  \emph{AIAA AVIATION 2022 Forum}, 2022, p. 3485.
\newblock \doi{10.2514/6.2022-3485}.

\bibitem[{Pitt(1980)}]{pitt1980rotor}
Pitt, D.~M., \emph{Rotor dynamic inflow derivatives and time constants from
  various inflow models}, Washington University in St. Louis, 1980.

\bibitem[{Pitt(1981)}]{pitt1981theoretical}
Pitt, D.~M., \enquote{Theoretical prediction of dynamic inflow derivatives,}
  \emph{Vertica}, Vol.~5, 1981, pp. 21--34.

\bibitem[{Hwang and Martins(2018{\natexlab{b}})}]{hwang2018fast}
Hwang, J.~T., and Martins, J.~R., \enquote{A fast-prediction surrogate model
  for large datasets,} \emph{Aerospace Science and Technology}, Vol.~75,
  2018{\natexlab{b}}, pp. 74--87.
\newblock \doi{10.1016/j.ast.2017.12.030}.

\bibitem[{Bouhlel et~al.(2019)Bouhlel, Hwang, Bartoli, Lafage, Morlier, and
  Martins}]{SMT2019}
Bouhlel, M.~A., Hwang, J.~T., Bartoli, N., Lafage, R., Morlier, J., and
  Martins, J. R. R.~A., \enquote{A Python surrogate modeling framework with
  derivatives,} \emph{Advances in Engineering Software}, 2019, p. 102662.
\newblock \doi{https://doi.org/10.1016/j.advengsoft.2019.03.005}.

\bibitem[{Drela(1989)}]{drela1989xfoil}
Drela, M., \enquote{XFOIL: An analysis and design system for low Reynolds
  number airfoils,} \emph{Low Reynolds number aerodynamics}, 1989, pp. 1--12.

\bibitem[{Selig(2014)}]{selig2014real}
Selig, M.~S., \enquote{Real-time flight simulation of highly maneuverable
  unmanned aerial vehicles,} \emph{Journal of Aircraft}, Vol.~51, No.~6, 2014,
  pp. 1705--1725.
\newblock \doi{10.2514/1.C032370}.

\bibitem[{Bryson et~al.(2016)Bryson, Marks, Miller, and
  Rumpfkeil}]{bryson2016multidisciplinary}
Bryson, D.~E., Marks, C.~R., Miller, R.~M., and Rumpfkeil, M.~P.,
  \enquote{Multidisciplinary design optimization of quiet, hybrid-electric
  small unmanned aerial systems,} \emph{Journal of Aircraft}, Vol.~53, No.~6,
  2016, pp. 1959--1963.
\newblock \doi{10.2514/1.C033455}.

\bibitem[{Lee et~al.(October, 2021)Lee, Ayton, Bertagnolio, Moreau, Chong, and
  Joseph}]{Lee:2021:PAS}
Lee, S., Ayton, L., Bertagnolio, F., Moreau, S., Chong, T.~P., and Joseph, P.,
  \enquote{Turbulent Boundary Layer Trailing-Edge Noise: Theory, Computation,
  Experiment, and Application,} \emph{Progress in Aerospace Sciences}, Vol.
  126, October, 2021, p. 100737.
\newblock \doi{10.1016/j.paerosci.2021.100737}.

\bibitem[{Li and Lee(October, 2020)}]{Li:2020:AHS}
Li, S., and Lee, S., \enquote{Prediction of Rotorcraft Broadband Trailing-Edge
  Noise and Parameter Sensitivity Study,} \emph{Journal of the American
  Helicopter Society}, Vol.~65, No.~4, October, 2020, pp. 1--14.
\newblock \doi{10.4050/JAHS.65.042006}.

\bibitem[{Li and Lee(July, 2021)}]{Li:2021:JAHS}
Li, S., and Lee, S., \enquote{Prediction of Urban Air Mobility Multi-Rotor VTOL
  Broadband Noise Using UCD-QuietFly,} \emph{Journal of the American Helicopter
  Society}, Vol.~66, No.~3, July, 2021, p. 032004.
\newblock \doi{10.4050/JAHS.66.032004}.

\bibitem[{Li and Lee(July, 2022)}]{Li:2022:JAHS}
Li, S., and Lee, S., \enquote{Acoustic Analysis and Sound Quality Assessment of
  a Quiet Helicopter for Air Taxi Operations,} \emph{Journal of the American
  Helicopter Society}, Vol.~67, No.~3, July, 2022, p. 032001.
\newblock \doi{10.4050/JAHS.67.032001}.

\bibitem[{Li and Lee(2022{\natexlab{a}})}]{Li:2022:AHS}
Li, S., and Lee, S., \enquote{Predictions and Validations of Small-Scale Rotor
  Noise Using UCD-QuietFly,} \emph{VFS}, 2022{\natexlab{a}}.

\bibitem[{Li and Lee(2022{\natexlab{b}})}]{Li:2022:AIAA}
Li, S., and Lee, S., \enquote{Extensions and Applications of Lyu and Ayton's
  Serrated Trailing-Edge Noise Model to Rotorcraft,} \emph{AIAA/CEAS},
  2022{\natexlab{b}}.

\bibitem[{Greenwood et~al.(2022)Greenwood, Brentner, Rau, and
  Gan}]{Greenwood:2022:IJA}
Greenwood, E., Brentner, K., Rau, R., and Gan, Z., \enquote{Challenges and
  Opportunities for Low Noise Electric Aircraft,} \emph{International Journal
  of Aeroacoustics}, Vol.~21, 2022, pp. 315--381.
\newblock \doi{10.1177/1475472X221107377}.

\bibitem[{Schlegel et~al.(1966)Schlegel, King, and
  Mull}]{schlegel1966helicopter}
Schlegel, R., King, R., and Mull, H., \enquote{Helicopter rotor noise
  generation and propagation,} Tech. rep., United Technologies Corp, Sikorsky
  Aircraft Division, 1966.

\bibitem[{Davidson and Hargest(1965)}]{Davidson:1965:RSS}
Davidson, I.~M., and Hargest, T.~J., \enquote{Helicopter Noise,} \emph{Journal
  of the Royal Aeronautical Society}, Vol.~69, 1965, pp. 325--336.
\newblock \doi{10.1017/s0001924000059583}.

\bibitem[{Johnson(2013)}]{johnson2013rotorcraft}
Johnson, W., \emph{Rotorcraft aeromechanics}, Vol.~36, Cambridge University
  Press, 2013.

\bibitem[{Li and Lee(2020)}]{li2020prediction}
Li, S.~K., and Lee, S., \enquote{Prediction of rotorcraft broadband
  trailing-edge noise and parameter sensitivity study,} \emph{Journal of the
  American Helicopter Society}, Vol.~65, No.~4, 2020, pp. 1--14.
\newblock \doi{10.4050/JAHS.65.042006}.

\bibitem[{Jia and Lee(2022)}]{jia2022computational}
Jia, Z.~H., and Lee, S., \enquote{Computational Study on Noise of Urban Air
  Mobility Quadrotor Aircraft,} \emph{Journal of the American Helicopter
  Society}, Vol.~67, No.~1, 2022, pp. 1--15.
\newblock \doi{10.4050/JAHS.67.012009}.

\bibitem[{Jia and Lee(April, 2022)}]{Jia:2022:JAHS}
Jia, Z., and Lee, S., \enquote{High-Fidelity Computational Analysis on the
  Noise of a Side-by-Side Hybrid VTOL Aircraft,} \emph{Journal of the American
  Helicopter Society}, Vol.~67, No.~2, April, 2022, p. 022005.
\newblock \doi{10.4050/JAHS.67.022005}.

\bibitem[{Sagaga and Lee(May 10--14, 2021)}]{Sagaga:2021:VFS}
Sagaga, J., and Lee, S., \enquote{Acoustic Predictions for the Side-by-Side Air
  Taxi Rotor in Hover,} \emph{VFS-77}, May 10--14, 2021.

\bibitem[{Gill et~al.(January 23--27, 2023)Gill, Lee, Ruh, and
  Hwang}]{Gill:2023:AIAA}
Gill, H., Lee, S., Ruh, M.~L., and Hwang, J.~T., \enquote{Applicability of
  Low-Fidelity Tonal and Broadband Noise Models on Small-Scaled Rotors,}
  \emph{AIAA}, January 23--27, 2023.

\bibitem[{Htwe(2019)}]{htwe2019design}
Htwe, N. H.~H., \enquote{Design of 50 kw Permanent Magnet Synchronous Motor for
  HEV,} \emph{IRE}, Vol.~3, 2019.

\bibitem[{Lu and Kar(2010)}]{lu2010review}
Lu, D., and Kar, N.~C., \enquote{A review of flux-weakening control in
  permanent magnet synchronous machines,} \emph{2010 IEEE Vehicle Power and
  Propulsion Conference}, IEEE, 2010, pp. 1--6.

\bibitem[{Zhao et~al.(2015)Zhao, Huang, Fang, and Li}]{zhao2015control}
Zhao, S., Huang, X., Fang, Y., and Li, J., \enquote{A control scheme for a High
  Speed Railway traction system based on high power PMSM,} \emph{2015 6th
  International Conference on Power Electronics Systems and Applications
  (PESA)}, IEEE, 2015, pp. 1--8.

\bibitem[{Gill et~al.(2005)Gill, Murray, and Saunders}]{gill2005snopt}
Gill, P.~E., Murray, W., and Saunders, M.~A., \enquote{SNOPT: An SQP algorithm
  for large-scale constrained optimization,} \emph{SIAM review}, Vol.~47,
  No.~1, 2005, pp. 99--131.
\newblock \doi{10.1137/S0036144504446096}.

\bibitem[{Gill et~al.(1981)Gill, Murray, and Wright}]{gills1981practical}
Gill, P.~E., Murray, W., and Wright, M., \emph{Practical Optimization. Ist
  Edn}, Academic Press, London, 1981.

\bibitem[{Patterson and Rao(2014)}]{patterson2014gpops}
Patterson, M.~A., and Rao, A.~V., \enquote{GPOPS-II: A MATLAB software for
  solving multiple-phase optimal control problems using hp-adaptive Gaussian
  quadrature collocation methods and sparse nonlinear programming,} \emph{ACM
  Transactions on Mathematical Software (TOMS)}, Vol.~41, No.~1, 2014, pp.
  1--37.
\newblock \doi{10.1145/2558904}.

\bibitem[{Wakefield and Dubuque(2009)}]{wakefield2009exceeding}
Wakefield, I., and Dubuque, C., \enquote{Exceeding tire speed rating during
  takeoff,} \emph{Boeing AERO Quarterly}, 2009.

\bibitem[{Holden and Goel(2016)}]{holden2016fast}
Holden, J., and Goel, N., \enquote{Fast-forwarding to a future of on-demand
  urban air transportation,} \emph{White Paper}, 2016.

\end{thebibliography}

\end{document}